\documentclass[aip,amsmath,amssymb,
preprint,%
]{revtex4-2}

\usepackage{graphicx}
\usepackage{dcolumn}
\usepackage{bm}
\usepackage{bbm}
\usepackage{subcaption} 
\usepackage{tikz}
\usepackage{standalone}
\usepackage{comment}

\newtheorem{lemma}{Lemma}

\usepackage[utf8]{inputenc}
\usepackage[T1]{fontenc}
\usepackage{mathptmx}
\usepackage{etoolbox}

\makeatletter
\def\@email#1#2{%
 \endgroup
 \patchcmd{\titleblock@produce}
  {\frontmatter@RRAPformat}
  {\frontmatter@RRAPformat{\produce@RRAP{*#1\href{mailto:#2}{#2}}}\frontmatter@RRAPformat}
  {}{}
}%

\renewcommand{\d}{{\rm d}}
\newcommand{\e}{{\rm e}}
\renewcommand{\i}{{\rm i}}

\newcommand{\FD}[2]{\frac{\d #1}{\d #2}}

\DeclareMathSymbol{\ZSet}{\mathalpha}{AMSb}{"5A}
\DeclareMathSymbol{\RSet}{\mathalpha}{AMSb}{"52}
\DeclareMathSymbol{\CSet}{\mathalpha}{AMSb}{"43}

\makeatother

\begin{document}

\title{Phase oscillator networks with multiple and state-dependent delays: 
A framework for exploring white matter plasticity in neurodynamics}

\author{G. Jolly}
\email[Corresponding author:]{grace.jolly@nottingham.ac.uk}
\author{R. Nicks}
\author{S. Ruschel}

\affiliation{School of Mathematical Sciences, University of Nottingham, Nottingham, NG7 2RD, UK}

\author{G. Iskenderoglu}
\affiliation{Department of Mathematics, Istanbul Ticaret University, Istanbul, 34840, T\"urkiye}

\author{S. Coombes}
\affiliation{School of Mathematical Sciences, University of Nottingham, Nottingham, NG7 2RD, UK}

\date{8 September 2026}

\begin{abstract}

Network science is increasingly focused on how node dynamics influence emergent phenomena such as oscillations, waves, chimeras, and turbulence. In oscillatory systems modeled by networks of coupled ordinary differential equations (ODEs), a common approach is to reduce the system to phase variables. When understanding how network delays shape emergent properties, the delays are often absorbed in the reduced description as a phase shift. The reduced system is a set of ODEs that loses some information about the full delay differential equation (DDE) system. We adopt a less restrictive approach and consider limit cycle oscillator networks with multiple delays that can be reduced to a DDE system. This captures the effects of delayed couplings, including coexistence and multistability of states. We analyze relative equilibria in the form of phase-locked states with tools previously developed for more general DDE settings. This allows us to explore patterning in networks with space-dependent delays, including in neuroscience, using symmetric bifurcation theory, linear stability analysis, and numerical simulations and continuation. In brain dynamics, time delays are determined by the speed of a communicating signal (action potential) along a fiber (axon). Importantly, these are now known to be state-dependent since the myelin (white matter) that insulates axons is plastic and can change in response to neuronal activity. A simple phenomenological model of this process (in the phase reduced description) is introduced and analyzed by extending techniques developed for fixed delays. Our analysis suggests that white matter plasticity can drive networks to more coherent behavior.

\end{abstract}

\maketitle

\begin{quotation}
The white matter that forms the brain axonal communication network is just as crucial to neurodynamics as the gray matter it interconnects. Increasingly, research shows that white matter is not merely a static substrate defining a fixed space-time connectome of linkages and delays; rather, these delays are mutable and subject to modulation. Transmission delays in neural communication arise from the finite speed at which action potentials propagate along axonal fibers a speed governed largely by the insulating properties of myelin, a white fatty substance. Crucially, myelin can be modified in response to neuronal activity, a phenomenon known as white matter plasticity.  This recognition compels the development of mathematical models and methodologies capable of handling brain-like networks with multiple, dynamically evolving delays. As a step in this direction, we investigate networks of coupled limit-cycle oscillators with delayed interactions, which can be reduced to lower-dimensional descriptions in terms of phase variables. While these reduced systems are still governed by delay differential equations, we show that analyzing oscillatory network states viewed as relative equilibria is significantly more tractable. Furthermore, our framework accommodates a biologically inspired form of white matter plasticity, wherein conduction delays become state-dependent. Exploring the consequences of this, we observe model behavior that aligns with growing empirical evidence: activity-dependent myelination can dynamically tune conduction speeds to enhance synchrony within neural circuits.

\end{quotation}

\section{Introduction\label{sec:Introduction}}

The study of large scale brain network models with delays is still in its infancy, at least with respect to theoretical investigations, although those based upon computational approaches are progressing rapidly, as exemplified by \cite{Rabuffo2021}. Typically this involves simulating a set of biologically inspired nodes arranged in a network with signals passed between nodes subject to a delay.  These nodes are often prescribed by oscillatory neural mass models, themselves constructed as possibly high-dimensional systems of ordinary differential equations, see e.g., \cite[Ch. 8]{Coombes2023}.  Readily available data such as that from the Human Connectome Project \cite{VanEssen2013} provides the strength of connection between nodes as well as the corresponding white matter tract lengths. When combined with information about the distribution of action potential speeds the latter can be used to estimate the time taken for signals to travel along axonal bundles.  Clearly, for this \textit{space-time} connectome there are as many delays in the network as there are connections.  Even for the simple modeling assumption of a common axonal propagation speed along the connecting fibers based on the 3D Euclidean distance between any 2 connected nodes in a 38 node Wilson--Cowan neural mass network, Deco \textit{et al}. have shown that white matter can play a key role in generating patterns of functional connectivity seen in so-called resting-state using structural data from the macaque brain \cite{Deco2009}.  Similar conclusions have been drawn when using human connectome data (with 66 anatomical nodes), although for a simplified phase-oscillator Kuramoto network \cite{Cabral2011,Budzinski2023,Thanos2023}.  
Moreover, Kuramoto networks with distance dependent delays have been shown to have explanatory power in understanding both cortical traveling waves \cite{Budzinski2023,Koller2024} and the frequency-dependent organization of functional networks \cite{Abolfazl2020}.
Such studies highlight the rich dynamics that can be exhibited by brain inspired network models with delays that are deeply relevant to modern theoretical neuroscience but remain only partially understood from a mathematical standpoint. From a neuroscience perspective this includes understanding how delays can control neuronal synchronization, see e.g., \cite{Protachevicz2020,Kusch2025}, and contribute to the
\textit{communication through coherence} hypothesis for cognition \cite{Fries2015}.  However, the challenge of how to understand the dynamics of large systems of interacting oscillators with many delays is not just specific to neuroscience and is relevant to physical models of lasers \cite{Soriano2013}, connected vehicle systems \cite{Szalai2013}, large ecosystems \cite{Pigani2022} and bio-inspired computation \cite{Esfahani2016,Tavakoli2024}, to name but a few other areas.  

One natural way to progress this challenge is to use the theory of weakly coupled oscillators \cite{Hoppensteadt97}, with a reduction to a phase description as in \cite{Ko2007,Ton2014}.  This often goes hand-in-hand with an approximation whereby delays manifest as a simple shift in a phase variable or as a shifted argument to a phase-interaction function, and see e.g., \cite{Ermentrout2009}.  The resulting system of ordinary differential equations (ODEs) avoids the analytical difficulty of dealing with infinite dimensional delay differential equations (DDEs), and a recent approach has been developed to improve over this first order approximation though still using ODEs \cite{Bick2024,Bick2025}.  Other (mean-field) approaches for treating phase-oscillator networks with delays have also been developed using the Ott--Antonsen (OA) ansatz \cite{Petkoski2016}, though this is only applicable when the phase interaction function is sinusoidal \cite{Petkoski2016}.  In this case, a further assumption of short delays can give rise to more tractable models, albeit with effective linear interactions only \cite{Zanette00,Perez2011}.  Here, we advocate for the use of phase reduced descriptions in coupled oscillator networks with delays, though not at the expense of forsaking a DDE description.  Although the analysis of periodic solutions of the full unreduced network equations is very challenging, with progress mainly through the use of numerical analysis, and see e.g., \cite{Coombes2025} for a recent harmonic balance approach, such solutions become \textit{relative equilibria} in the phase reduced framework. These can be readily constructed as phase-locked states and their stability can be analyzed with techniques for the study of standard equilibria in delayed systems.  This is a main topic for the work considered in this paper.  Moreover, we recognize the fact that delays in brain networks are subject to variation on long time scales due to the fact that conduction speeds depend on myelination of the axonal connections and these in turn adapt in response to neural activity.  With the introduction of a biologically inspired white matter plasticity rule we further use the framework presented to explore the consequences of such white matter plasticity on emergent network states.

In \S \ref{model} we introduce networks of weakly coupled oscillators as a convenient class of models of coupled oscillator networks with delays amenable to a phase description. This sets the scene for the phase oscillator network DDE model that we analyze in subsequent sections.   The theoretical work necessary for describing phase-locked states is presented in \S \ref{phaselocked} together with the machinery for determining their linear stability.  This is illustrated with bifurcation studies of a simple two-node network and larger and more general networks subject to a row-sum constraint on the matrix representation of the weighted graph of interactions.  
This general approach is used to treat the special case of ring networks with distant dependent delays in 
\S \ref{ring} using tools from symmetric bifurcation theory, and we also show how to treat twisted states.  A phenomenological rule for white matter plasticity is introduced in \S \ref{plasticity}, and the extension of the theoretical work from fixed to state-dependent delays is presented.  As well as exploring how this plasticity rule changes the emergent dynamics on idealized ring networks, we also consider networks built from human connectome data.  In both cases, using theory and numerical simulations, we find that white matter plasticity can drive networks to more coherent behavior. 
Finally in \S \ref{discussion} we summarize the main results and conclusions of the work presented in this paper and discuss natural extensions and next steps.

\section{Phase oscillator network model\label{model}}

We consider the dynamics of networks of $N$ identical oscillators $x_i$ with $1\leq i \leq N,$ $m \geq 2,$ and ${x}_i(t) \in \RSet^m$ given by
\begin{equation}
	\dot{x}_i (t)\equiv \FD{}{t} {x}_i (t)= {f}({x}_i(t)) + \sigma \sum_{j=1}^N  W_{ij}  {G} ({x}_i(t),{x}_j(t -\tau_{ij})), 
	\label{Network1}
\end{equation}
with pairwise interactions.

Here, $G({x}_i,{x}_j)$ is the dynamics that describes the coupling between nodes $i$ and $j$, the relative strength of this interaction is $W_{ij} \in \RSet$,
and $\sigma$ sets the overall network coupling strength.  Importantly, we account for delays in the interactions between nodes with the inclusion of lags in the interaction defined by $\tau_{ij} \geq 0$.  Thus, the system of equations (\ref{Network1}) is a set of DDEs with a heterogeneous space-time connectome.
The space-time connectome specifying the interaction between nodes $i$ and $j$ decomposes into two parts:  i) the structural connectome specified by the weights $W_{ij}$, and ii) the corresponding temporal delays specified by $\tau_{ij}$.  
We may think of each oscillator as associated with a node of a structural network (which, most traditionally, takes the form of a graph \cite{Newman2010}), and each delayed interaction is associated with an edge of that network.
For a system of uncoupled nodes ($\sigma=0$) we shall assume that the local node dynamics described by the ODE $\dot{s} = f(s)$ admits a $T$-periodic attracting hyperbolic limit cycle.  In this case, the theory of weakly-coupled oscillators, see e.g.,  Ref. \onlinecite{Hoppensteadt97}, is expected to give rise to a reduced description from $Nm$ DDEs to $N$ DDEs of the form
\begin{equation}
\FD{}{t} \theta_i(t) = \omega + \sigma \sum_{j=1}^N W_{ij} H(\theta_{j}(t-\tau_{ij}) - \theta_i(t)), 
\label{PhaseNetwork}
\end{equation}
with $\theta_i(t) \in [0, 2 \pi)$ and where $\omega=2\pi/T$ is the frequency of the $i$th node when uncoupled from the rest of the network. For completeness we include a derivation of (\ref{PhaseNetwork}) from (\ref{Network1})  in Appendix \ref{phasereduction}.
In (\ref{PhaseNetwork}), $H$ is a $2 \pi$ periodic \textit{phase interaction function} given by (\ref{H}) in Appendix \ref{phasereduction}.  
However, for the purposes of this paper we shall fixate on the illustrative choice of a bi-harmonic function \cite{Ashwin2016}:
\begin{equation}
H(\theta) = -\sin(\theta - a) + r \sin(2 \theta) ,
\label{biharmonic}
\end{equation}
and treat the pair $(a,r)$ as bifurcation parameters.  For the choice $(a,r)=(\pi,0)$ we recover the Kuramoto model.

It is worthwhile to note that the phase-shift approximation, whereby $\theta_i(t-\tau) \simeq \theta_i(t) - \omega \tau$, yields the ODE system
\begin{equation}
\dot{\theta}_i = \omega + \sigma \sum_{j=1}^N W_{ij} H(\theta_{j} - \theta_i - \omega \tau_{ij}), 
\label{PhaseShift}
\end{equation}
which is expected to lose some information about network behavior compared to the DDE systems given by (\ref{PhaseNetwork}).  For example, time delays have been reported to induce multistability of synchronized states in the Kuramoto model, a feature which is not retained under phase-shift approximation \cite{yeung1999time}. It is also important to note that the number of delays can potentially be reduced by imposing a suitable component-wise time-shift transformation \cite{lucken2013reduction}.

\section{Phase-locked states\label{phaselocked}}

Here we consider states in which all pairs of oscillators are frequency-locked with a constant phase lag between each pair.
We denote such a phase-locked state by $\theta_{i}(t) = \Omega t + \phi_i$ for some \textit{emergent frequency} $\Omega$ and time-independent phases $\phi_i$.
Fixing $\phi_1=0$ as a reference phase, the set of $N$ unknowns $(\Omega, \phi_2, \ldots, \phi_N)$ is found as the solution to the system of $N$ nonlinear equations:
\begin{equation}
\Omega = \omega + \sigma \sum_{j=1}^N W_{ij} H(\phi_j - \phi_i -\Omega \tau_{ij}) .
\label{emergent}
\end{equation}
Given the nonlinear nature of this system of equations it is pragmatic to use numerical routines based on Newton--Krylov methods when constructing branches of solutions.  Note that the phase-shift approximation (\ref{PhaseShift}) would yield an explicit emergent frequency $\Omega = \omega + \sigma \sum_{j=1}^N W_{ij} H(\phi_j - \phi_i -\omega \tau_{ij})$,  without the possibility of multiple solutions for $\Omega$ for a fixed set of phases.

\subsection{Linear stability\label{sec:linstab}}

Here, we treat the linear stability of phase-locked states by introducing small perturbations $\delta \theta_i(t)$ and linearizing (\ref{PhaseNetwork}) by writing $\theta_{i}(t) = \Omega t +\phi_i + \delta \theta_i(t)$.  This gives
\begin{equation}
\FD{}{t} \delta \theta_i(t) = \sigma \sum_{j=1}^N W_{ij} H'(A_{ij}) \left [ \delta \theta_{j}(t-\tau_{ij}) - \delta \theta_i(t) \right ], 
\end{equation}
where $A_{ij} = \phi_j - \phi_i -\Omega \tau_{ij}$.
For separable solutions of the form $\delta \theta_i(t) = u_i \e^{\lambda t}$, $\lambda \in \CSet$, this yields
\begin{equation}
\lambda u_i =  \sum_{j=1}^N  \left [ \mathcal{J}_{ij}(\lambda) u_{j} - \mathcal{J}_{ij}(0) u_i \right ] , 
\label{ui}
\end{equation}
where
\begin{equation}
 \mathcal{J}_{ij}(\lambda) = \sigma W_{ij} H'(A_{ij}) \e^{-\lambda \tau_{ij}} .
 \label{J}
\end{equation}
Thus, solutions will be stable provided $\text{Re} \, (\lambda) <0$ where $\lambda$ is a zero of the function $\mathcal{E}(\lambda) = \det [\lambda I_N + \mathcal{L}(\lambda)]$, where $\mathcal{L}(\lambda)$ has a graph-Laplacian like structure with
\begin{equation}
\mathcal{L}_{ij}(\lambda) = - \mathcal{J}_{ij} (\lambda) + \delta_{ij} \sum_k \mathcal{J}_{ik}(0) .
\end{equation}
Here $I_N$ is the $N \times N$ identity matrix.  
We note that the $N$-dimensional vector $(1,1,\ldots,1)$ is an eigenvector of $\mathcal{L}(0)$ with eigenvalue $0$.  Hence, 
$\lambda=0$ is a always a solution of $\mathcal{E}(\lambda) = 0$.  This reflects the invariance of (\ref{emergent}) to phase-shifts.

In general for $\lambda \neq 0$ the matrices $\mathcal{J} (\lambda)$ and $\mathcal{J} (0)$ will not have a common eigenspace.  However, this is the case for a single delay ($\tau_{ij} = \tau$ for all $(i,j)$) and for systems with space-dependent interactions on a ring so that the matrices with components $W_{ij}$ and $\tau_{ij}$ are circulant (assuming solutions where $H'(A_{ij})$ is also circulant).  For these cases, after introducing right (left) normalized eigenvectors of $\mathcal{J}(0)$ as $a^\mu$ ($b^\mu$), we can write $\mathcal{J}(\lambda) = \sum_{\mu} \gamma_\mu(\lambda) b^\mu \otimes a^\mu$, where $\otimes$ denotes the tensor product.  Here,  $\gamma_\mu(\lambda)$ can be constructed by projection of $\mathcal{J}(\lambda)$ on to $a^\mu$ and $b^\mu$ as
$\gamma_\mu(\lambda) = \left ( a^\mu \right )^\top \mathcal{J}(\lambda) b =  \sum_{i,j} a^\mu_i \mathcal{J}_{ij}(\lambda) b^\mu_j$.  
By introducing the set of row-sums $k_i = \sum_{j=1}^N \mathcal{J}_{ij}(0)$ and considering the vector $u=(u_1,\ldots,u_N)$ to be parallel to a right eigenvector of $\mathcal{J}(0)$ we can diagonalize (\ref{ui}) to the form
\begin{equation}
(\lambda + k_i - \gamma_\mu(\lambda)) u_i =  0, \qquad \mu =1,\ldots, N .
\label{uidiag}
\end{equation}
For non-trivial solutions of (\ref{uidiag}) we must have that $\mathcal{E}_\mu(\lambda) = 0$, where
\begin{equation}
\mathcal{E}_\mu(\lambda) = \det \left [
(\lambda - \gamma_\mu(\lambda)) I_N + \operatorname{diag}(k_1, \ldots, k_N)
\right ] .
\label{E}
\end{equation}
Equivalently, we may write $\mathcal{E}_\mu(\lambda) = \prod_{i=1}^N (\lambda-\gamma_\mu(\lambda)+k_i)$.
In this way we may determine the stability of a phase locked state according to the condition $\text{Re}\, (\lambda) <0 $ for all $\mu$ (excluding the zero eigenvalue that arises from phase-shift symmetry).

\subsection{A simple two oscillator network\label{N=2}}

It is instructive to first consider the case of two identical oscillators with reciprocal and symmetric coupling. In this case we set
$\tau_{12}  = \tau = \tau_{21}$,  and
\begin{equation}
W = \begin{bmatrix} 0 & 1 \\
1 & 0
\end{bmatrix} .
\end{equation}
Now consider the phase-locked state with phase difference $\psi = \phi_2  - \phi_1$. The emergent frequency is given by the equations
\begin{align}
    \Omega &= \omega + \sigma H(\psi - \Omega \tau), \nonumber \\
    \Omega &= \omega + \sigma H( -\psi - \Omega \tau).
\label{eq: emergent2_all}
\end{align}
It is worthwhile to note that because of the periodicity of $H$ the system (\ref{eq: emergent2_all}) defining a solution branch $\Omega = \Omega(\tau)$ is invariant under a shift in $\tau$ such that $\Omega \tau \mapsto \Omega \tau + 2 \pi \ZSet$.  Thus, if a solution pair $(\tau, \Omega)$ can be found then other solutions $(\tau+2 \pi/\Omega, \Omega)$ will be generated.  This is consistent with more general observations on the \emph{reappearance} of periodic solutions in delayed systems \cite{Yanchuk2009}.

We first consider the synchronous, $\psi=0$, and anti-synchronous, $\psi = \pi$, solutions.  In these cases the emergent frequency satisfies one equation given by 
\begin{equation}
\Omega = \omega + \sigma H(\psi - \Omega \tau) .
\label{emergent2}
\end{equation}
We use the prescription from Sec.~\ref{sec:linstab} to consider the stability of the synchronous and anti-synchronous solutions. We have that $\gamma_1(\lambda) = \beta \e^{-\lambda \tau} = - \gamma_2(\lambda)$ and $k_1 = \beta = k_2$, where
\begin{equation}
\beta = \sigma H'(\psi -\Omega \tau) .
\label{alpha}
\end{equation}
Using (\ref{E}), the eigenvalues that determine stability are determined by the two equations 
\begin{equation}
\mathcal{E}_\pm(\lambda) = \lambda + \beta \left [ 1 \pm \e^{-\lambda \tau} \right ] = 0 .
\label{Epm}
\end{equation}

For the synchronous and anti-synchronous solutions it is straightforward to show that there are no \textit{dynamic} instabilities.  By this we mean that eigenvalues cannot cross through the imaginary axis away from the real line.  To see this consider a search for solutions to the spectral problem defined by $\mathcal{E}_\pm(\i \varpi) = 0$ with $\varpi \in \RSet / 0$.  Substitution of $\lambda = \i \varpi$ into (\ref{Epm}) and balancing real and imaginary terms yields the pair of equations
\begin{align}
0 & = \beta \left [ 1 \pm \cos(\varpi \tau) \right ] , \label{i} \\
0 & = \varpi \mp \beta \sin (\varpi \tau)   . \label{ii}
\end{align}

For $\beta \neq 0$ (\ref{i}) has solutions when $\varpi \tau = \ZSet \pi$, and then using (\ref{ii}) this gives $\varpi = 0$.  Given that there are no dynamic instabilities we now restrict attention to $\lambda \in \RSet$.  In this case $\mathcal{E}_\pm (\lambda) \in \RSet$ and we note the following:  $\mathcal{E}_+ (0) = 2 \beta$, $\mathcal{E}_- (0) = 0$, 
$\mathcal{E}'_\pm (0) = 1 \mp \beta \tau$, and $\lim_{\lambda \rightarrow \infty} \mathcal{E}_{\pm}(\lambda) = \lambda$.  For $\beta >0$, $\mathcal{E}_+(\lambda) > \mathcal{E}_-(\lambda) >0$ for all $\lambda >0$. Hence, there are no solutions of $\mathcal{E}_\pm(\lambda) = 0$ with positive $\lambda$.  For $\beta <0$ the graph of $\mathcal{E}_+(\lambda)$ must cross the zero axis at a positive value of $\lambda$.  Hence, a solution (synchronous or anti-synchronous) is linearly stable if $\beta >0$, where $\beta$ is given by (\ref{alpha}).

It is also possible to establish that changes in stability along synchronous and anti-synchronous solution branches, $\Omega = \Omega(\tau)$, defined by (\ref{emergent2}) occur at stationary points.
To see this, note that differentiation of (\ref{emergent2}) with respect to $\tau$ gives
\begin{equation}
\FD{\Omega}{\tau} = - \frac{\beta \Omega}{1+\beta\tau} .
\end{equation}
Hence, a bifurcation defined by $\beta=0$ coincides with a stationary point of the curve $\Omega = \Omega(\tau)$.


\begin{figure}
    \begin{subfigure}[t]{0.6\columnwidth}
        \centering
        \includegraphics[width=\linewidth]{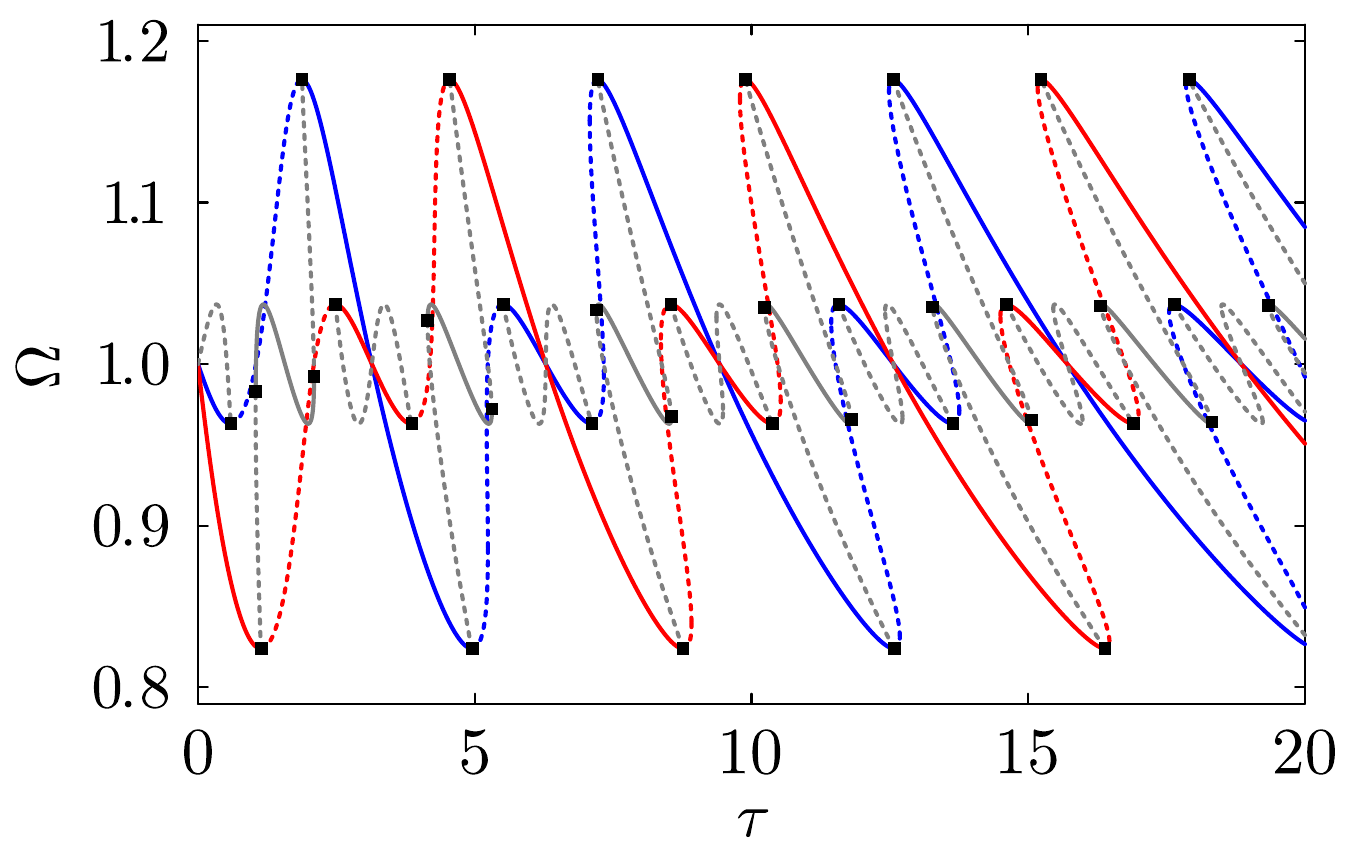}
        \caption{Emergent frequency $\Omega$ against delay $\tau$.}
        \label{subfig:N=2 a}
    \end{subfigure}

    \medskip

    \begin{subfigure}[t]{0.6\columnwidth}
        \centering
        \includegraphics[width=\linewidth]{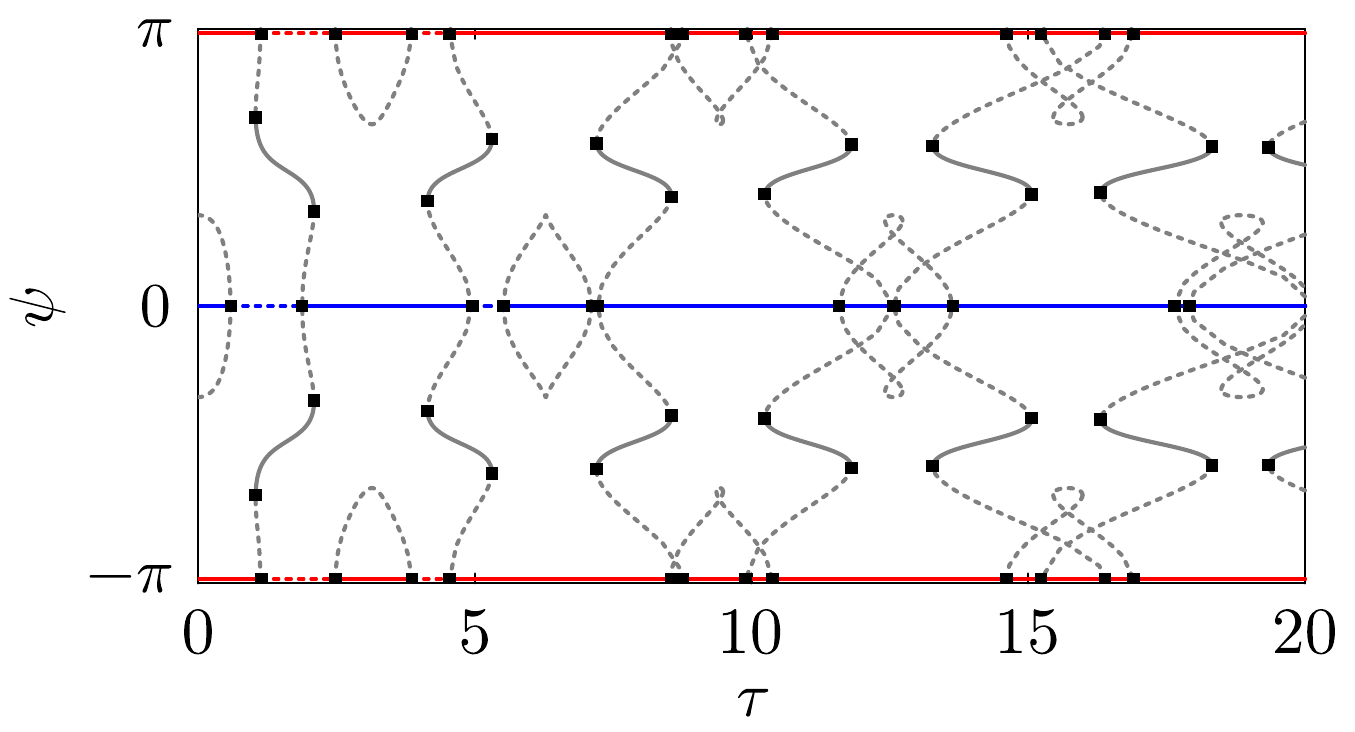}
        \caption{Phase difference $\psi$ against delay $\tau$.}
        \label{subfig:N=2 b}
    \end{subfigure}

    \caption{Emergent frequency $\Omega$ and phase difference $\psi$ as functions of delay $\tau$ in a two node network. Synchronous solutions ($\psi=0$) are rendered in blue, anti-synchronous solutions ($\psi=\pi$) in red and intermediate solutions ($\psi \ne 0, \pi$) in gray.
    Stable (unstable) solution branches are indicated by solid (dashed) lines.
    Black squares denote changes in stability associated with a real eigenvalue crossing zero. 
    Stability changes occur at stationary points of $\Omega=\Omega(\tau)$ along the synchronous and anti-synchronous branches shown in (a), and no dynamical instabilities are observed along these branches. Intermediate solution branches emerge from bifurcation points on the synchronous or anti-synchronous branches. 
    Parameters:  $\omega=1$, $a=0$, $r=1$, $\sigma=0.1$.}
    \label{Fig:N=2}

\end{figure}

We plot the solution branches $\Omega = \Omega(\tau)$ and $\psi = \psi(\tau)$ for all phase-locked solutions with $\phi_2  - \phi_1 = \psi$. We plot these curves parametrically and compute the stability numerically, as described in Appendix \ref{appendix parametric plotting two node}. An example is shown in 
Fig.~\ref{Fig:N=2}. 
Fig.~\ref{subfig:N=2 a} shows that for small $\tau$ there is a single branch of each type of solution.  However, for sufficiently large $\tau$ multiple solutions of the same type may stably co-exist.  For example, for $\tau=20$ in Fig.~\ref{subfig:N=2 a} three stable synchronous solutions co-exist (distinguished by their emergent frequency). No dynamic instabilities occur along the synchronous and anti-synchronous solution curves, which agrees with the linear stability analysis. Also, the stability of these solutions changes at the stationary points of $\Omega= \Omega(\tau)$ of the corresponding curves. Intermediate phase locked states given by $\psi \ne 0, \pi$ emerge from bifurcation points on the synchronous or anti-synchronous curves. Fig.~\ref{subfig:N=2 b} shows the possible values of the phase difference $\psi$ for each delay value. If $\psi$ satisfies (\ref{eq: emergent2_all}), then $-\psi$ also satisfies the system. This symmetry is illustrated by the reflection symmetry about $\psi = 0$ in Fig.~\ref{subfig:N=2 b}.

\begin{figure*}
    \centering

    \includegraphics[width=0.95\linewidth]{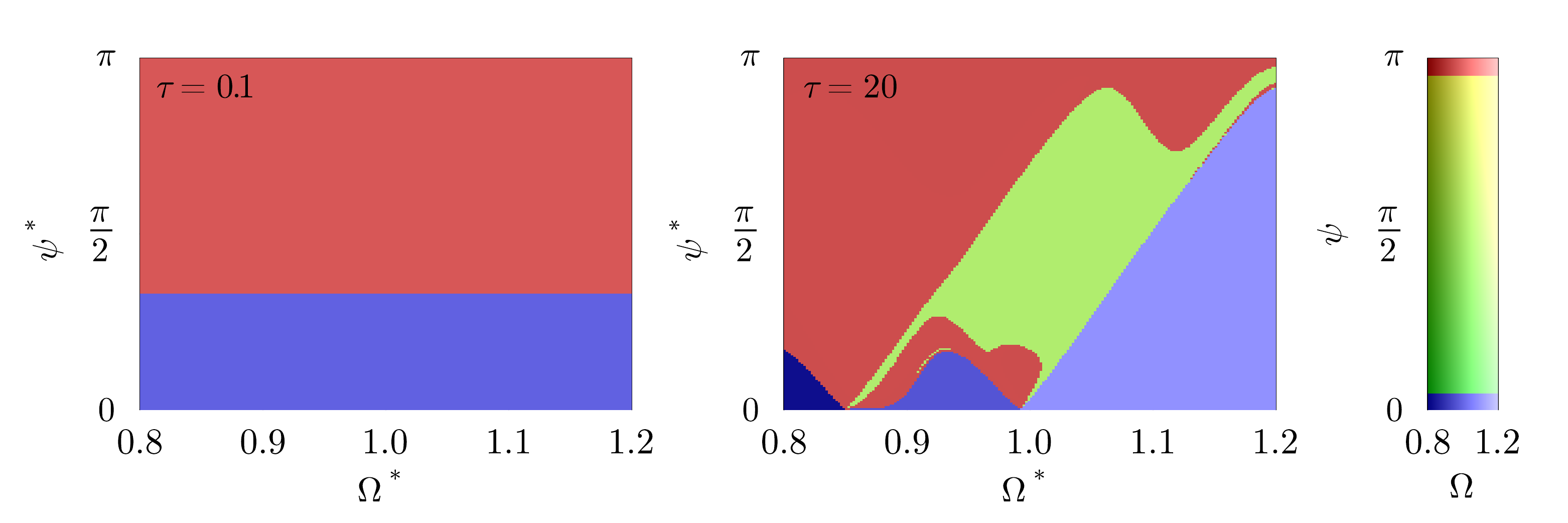}

    \caption{Heatmaps showing the eventual emergent frequency $\Omega$ and phase difference $\psi$ in a two node network, for varying initial conditions with initial frequency $\Omega^*$ and initial phase difference $\psi^*$. 
    Left: At a small delay $\tau = 0.1$. Right: At a large delay $\tau = 20$.
    Simulations are performed on a $200\times200$ grid up to $t=800$, with history functions $h_1(t)=\Omega^* t$ and $h_2(t)=\Omega^* t + \psi^*$. Other parameters as in Fig.~\ref{Fig:N=2}.}
    \label{Fig:N2_heatmap}
\end{figure*}

We investigate the basins of attraction of the two node network at two delay values in Fig.~\ref{Fig:N2_heatmap}, initializing the system with history functions $h_1(t)=\Omega^* t$ and $h_2(t)=\Omega^* t + \psi^*$ over a range of values of $(\Omega^*, \psi^*)$. The left panel shows that at a small delay the system tends to one of two states: initial phase differences close to synchrony ($\psi^* = 0$) lead to a synchronous solution, and those further away lead to an anti-synchronous solution. In both cases the emergent frequency is close to 1. The initial frequency does not affect the resulting state. This is consistent with Fig.~\ref{Fig:N=2} which shows that for delays close to zero there are two stable solutions, one synchronous and one anti-synchronous, both with emergent frequency close to 1.
The second plot shows that with a longer delay $\tau = 20$ the system can evolve to one of five different states: three synchronous solutions with different emergent frequencies shown by three shades of blue, one anti-synchronous solution, and one intermediate phase locked solution. This is also consistent with Fig.~\ref{Fig:N=2} which shows at $\tau=20 $ there are five stable states which match these types of solutions. In this case, both the initial frequency $\Omega^*$ and phase difference $\psi^*$ influence the resulting state.

\subsection{Synchrony in networks with a row-sum constraint and a single delay}

We consider now the case of an arbitrary number of identical nodes albeit with a single fixed delay, $\tau_{ij} = \tau$ for all $(i,j)$ and a row-sum-constraint $\sum_{j=1}^N W_{ij} = \Gamma$ for all $i$.
The governing equations for a synchronous solution $\theta_i(t) = \Omega t$ for all $i$ reduce to
\begin{equation}
\Omega = \omega + \sigma \Gamma H(-\Omega \tau).
\end{equation}
Moreover, the form for $k_i$ reduces to $k_i=k$ for all $i$, where $k = \sigma \Gamma H'(-\Omega \tau)$, and $\gamma_{\mu}(\lambda) = \sigma  H'(-\Omega \tau) \e^{-\lambda \tau} w_\mu$, where $w_\mu$ is an eigenvalue of $W$.  Hence from \eqref{E}, $\mathcal{E}_\mu (\lambda)=0$ when
\begin{equation}
\lambda + k \left [ 1 - \frac{w_\mu}{\Gamma} \e^{-\lambda \tau} \right ] = 0 .
\label{lambdarowsum}
\end{equation}
We note that the $N$-dimensional vector $(1,1,\ldots,1)$ is an eigenvector of $W$ with eigenvalue $\Gamma$, which ensures that $\lambda=0$ is always a solution of (\ref{lambdarowsum}) as expected.  The study above essentially recovers a previous analysis of Earl and Strogatz \cite{Earl2003} who established that synchrony is linearly stable if $k >0$.  
An identical argument to that used in Sec.~\ref{N=2} can be used to establish that a bifurcation defined by $k=0$ coincides with a stationary point of the curve $\Omega = \Omega(\tau)$.

Although seemingly restrictive, the choice of a row-sum constraint covers the case of globally coupled networks, which is often invoked in the study of Kuramoto networks since it allows a mean-field reduction using the OA ansatz.  Indeed there are now several studies of (\ref{PhaseNetwork}) for $W_{ij} = N^{-1}$ and $H(\theta)=\sin \theta$ incorporating a distribution of interaction delays (and not just a single one as considered in this section) in the thermodynamic limit $N \rightarrow \infty$ that further highlight the impact delays can have on Kuramoto networks, see e.g., \cite{Lee2009,Niu2014,Niu2017,Skardal2018}.  We shall not pursue this further here, and instead turn to more structured systems, albeit arranged on a ring.

\section{Ring networks and distant dependent interactions\label{ring}}

Here we consider a coupled oscillator system with interactions whose strength and delay both depend on distance.  For simplicity we consider nodes arranged on a ring with the distance between nodes $i$ and $j$ given by  $\operatorname{dist} (i,j) = \min(|i-j|, N-|i-j|)d$, for some spatial scale set by $d$.  The space-time connectome then has a circulant matrix representation with $W_{ij} = W_{\operatorname{dist} (i,j)}$ and $\tau_{ij} = \tau_{\operatorname{dist} (i,j)}$ with rows
generated by $W_{j} = w(\operatorname{dist} (0,j))$ for $j=0, \ldots,N-1$ for some function $w$ and 
\begin{equation}
\tau_{j}=
\begin{cases}
\frac{\operatorname{dist} (0,j)}{v} & j \neq 0 \\
\tau_0 & j=0
\end{cases}.
\label{tau_spacedependent}
\end{equation}
Here the delays $\tau_{j}$ for $j \neq 0$ are determined in terms of a single fixed conduction speed $v >0$.  The \textit{self-delay} $\tau_0 \geq 0$ is independent of $v$.

There is a synchronous solution with emergent frequency defined by
\begin{equation}
\Omega =\omega + \sigma \sum_{j=0}^{N-1} W_j H(-\Omega \tau_j)  .
\label{eq: synch ring}
\end{equation}
Another solution that is readily constructed for such an architecture is that of a splay solution with $\theta_i(t) = \Omega t + (2 \pi i)/N$, with emergent frequency
\begin{equation}
\Omega =\omega + \sigma  \sum_{j=0}^{N-1} W_j H(2 \pi j/N-\Omega \tau_j)  .
\label{omega_splay}
\end{equation}
In both cases (synchrony and splay) we note that a solution with emergent frequency $\Omega$ for a fixed value of $v$ will \textit{reappear} when $\Omega d/v \mapsto \Omega d /v + 2 \pi \ZSet$.

The constructs of \S \ref{sec:linstab} that help determine the linear stability similarly take more simplified circulant forms with  $A_{ij} = A_{\operatorname{dist} (i,j)}$, $\mathcal{J}_{ij}(\lambda) = \mathcal{J}_{\operatorname{dist} (i,j)}(\lambda)$, $k_i=k$ for all $i$ where $k = \sigma \sum_{j=0}^{N-1} W_j H'(A_{\operatorname{dist} (0,j)})$.  Most usefully, we may exploit the properties of circulant matrices to obtain
\begin{equation}
\gamma_\mu(\lambda) = \sum_{j=0}^{N-1} \mathcal{J}_{\operatorname{dist} (0,j)}(\lambda) \exp \left (
\frac{2 \pi {\rm i} (\mu-1) j}{N}
\right ) ,
\label{eq:stability splay}
\end{equation}
where $\mathcal{J}_{\operatorname{dist} (0,j)}(\lambda) = \sigma W_j H'(\psi_j -\Omega \tau_j)\e^{-\lambda \tau_j}$, where $\psi_i = 0$ for the synchronous solution and $\psi_i = 2 \pi i /N$ for the splay state.  The eigenvector associated with eigenvalue $\gamma_\mu$ is given by
$p_\mu = (1,a_\mu, a_\mu^2, \ldots, a_\mu^{N-1})/\sqrt{N}$, where $\mu=1,\ldots, N$, and $a_\mu=\exp(2 \pi \i (\mu -1)/N)$ are the $N$th roots of unity.

Given that there are changes in $v$ associated with aging \cite{Petkoski2023} and multiple-sclerosis \cite{Marti-Juan2023} it is interesting to explore the dynamical behavior of (\ref{PhaseNetwork}) with space-dependent delays prescribed by (\ref{tau_spacedependent}) treating $v$ as a primary bifurcation parameter. We plot solution curves in the form $\Omega = \Omega(v^{-1})$ to follow the structure of the frequency-delay diagram in Fig.~\ref{Fig:N=2} , since $v^{-1}$ is proportional to the delays in the ring network. We again plot these curves parametrically, using the method described in Appendix \ref{appendix parametric plotting ring}.

\begin{figure}
\centering
\includegraphics[width=0.55\linewidth]{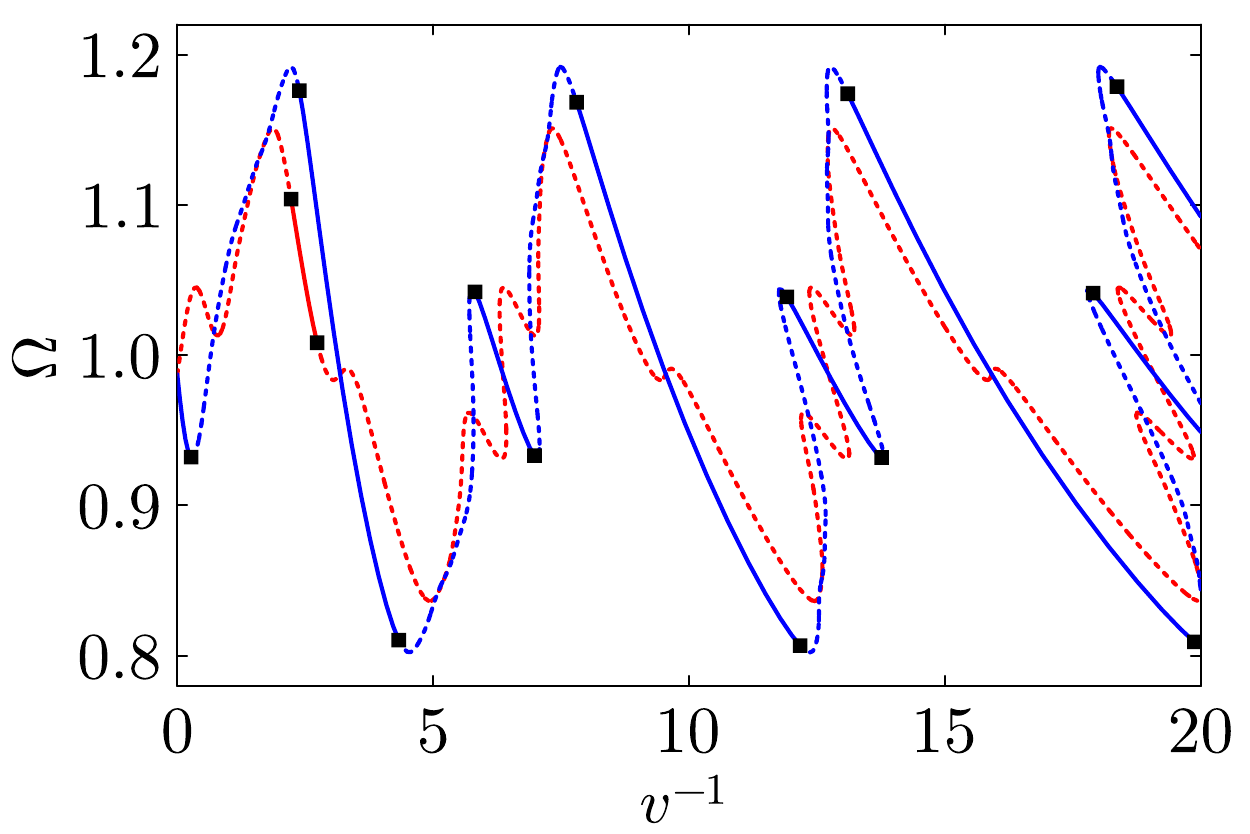}
\caption{Emergent frequency $\Omega$ for synchronous (blue) and splay (red) states as a function of $v^{-1}$ in a ring network with distant dependent interactions. Stable (unstable) solution branches are indicated by solid (dashed) lines.
Black squares denote changes in stability associated with a real eigenvalue crossing zero. Here $N=11$ and $w(x)=\exp(-|x|)/2$.
Parameters:  $\omega=1$, $a=0$, $r=1$, $\sigma=0.3$, $d=1$, $\tau_0=1$.
\label{Fig:RingSyncandSplay}
}
\end{figure}

Fig.~\ref{Fig:RingSyncandSplay} shows the solutions curves $\Omega = \Omega(v^{-1})$ for the synchronous and splay states in an $N=11$ node network. Other phase locked solutions exist for this network but are not shown on this diagram. Similarly to the two node network, for small delays (small values of $v^{-1}$), there is a single branch of each type of solution, but for larger values of $v^{-1}$ multiple solutions of the same type may exist. For example, at $v^{-1} = 19$ in Fig.~\ref{Fig:RingSyncandSplay} there are three stable synchronous solutions. There are also no dynamic instabilities along the branches, but the changes of stability do not occur at stationary points. There are values of $v^{-1}$ for which neither the synchronous or the splay solution is stable, for example at $v^{-1}=5$. This suggests the network does not evolve to synchrony or the splay state at these parameter values. It may evolve to other phase locked states, some of which are discussed in the remainder of this section.

\subsection{Twisted States}

One type of phase locked solution that the ring network admits is uniformly twisted states where there is a constant phase difference between consecutive oscillators. For $q = 0, \dots, N-1$, 
the $q$-twisted state is the phase-locked solution $\theta_i(t) = \Omega t + \phi_i^{(q)}$ with
\begin{equation}
    \phi_i^{(q)} = \frac{2\pi q i}{N},
\end{equation}
where $q$ is the winding number which measures the number of full twists in phase when you go around the ring once. The case $q=0$ corresponds to synchrony, and $q=1$ corresponds to the splay state discussed above.

Proceeding analogously to the construction of (\ref{omega_splay}) for a splay state, we find that the emergent frequency for a twisted state is
\begin{equation}
\Omega = \omega + \sigma \sum_{j=0}^{N-1} W_j 
H\!\left(\frac{2 \pi q j}{N} - \Omega \tau_j \right).
\label{omega_twisted}
\end{equation}
The linear stability is the same as in (\ref{eq:stability splay}) with $\psi_i = 2\pi q i/N$, and the solution branches $\Omega = \Omega(v^{-1})$ can be plotted parametrically following the method in Appendix \ref{appendix parametric plotting ring}.

\begin{figure}
    \centering
    \includegraphics[width=0.55\linewidth]{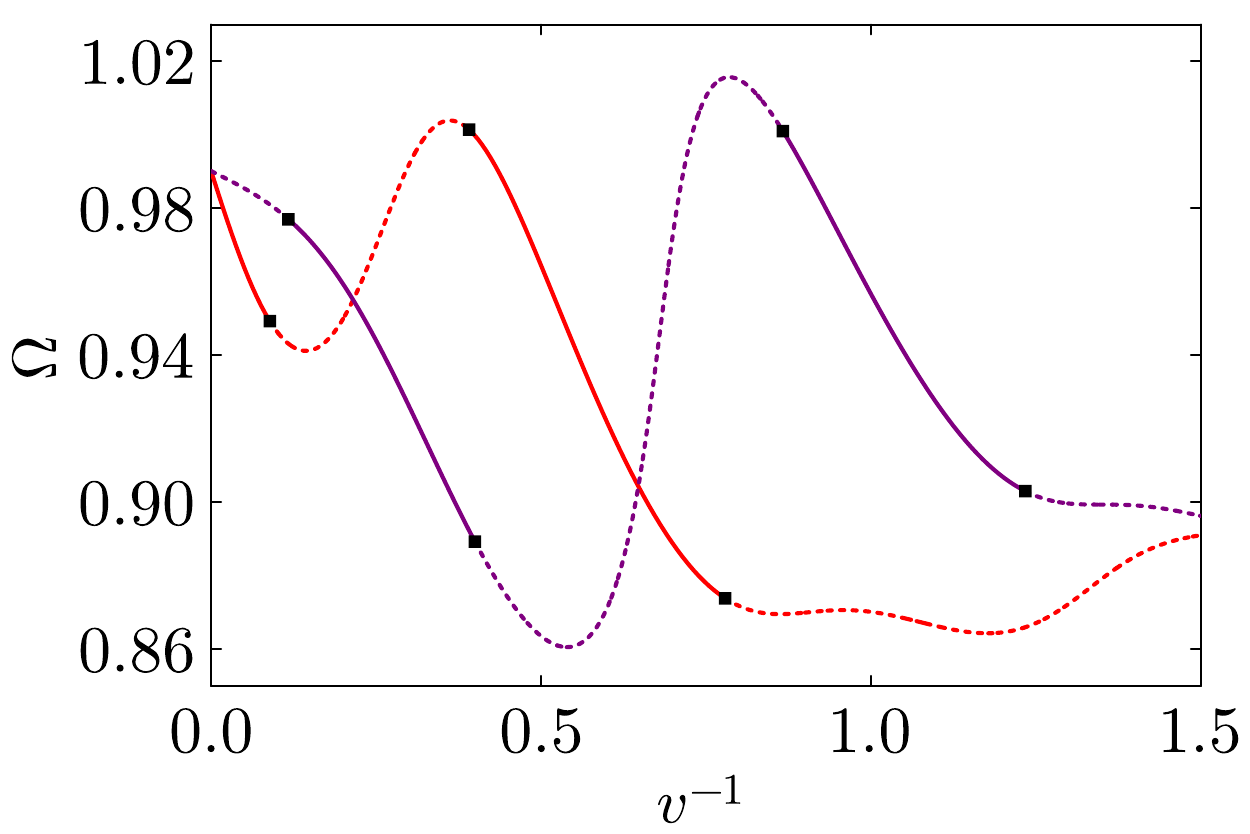}
    \caption{ Emergent frequency $\Omega$ as a function of $v^{-1}$ for twisted states in an $N=30$ network. The splay state is rendered in red and the $3$-twisted state in purple. Linearly stable solutions are indicated by solid lines and unstable solutions by dashed lines. Connectivity is given by the wizard hat function $w(x) = \left( 1 - |x|/3\right) \exp\left(-|x|/3\right)$. Parameters: $\omega=1$, $a= \pi$, $r=0$, $\sigma=0.1$, $d=1$, $\tau_0=0.1$.   For certain values of $v^{-1}$ the splay state is unstable but the $3$-twisted state is stable. }
    \label{fig:twisted branches}
\end{figure}

\begin{figure}
    \begin{subfigure}[t]{\columnwidth}
        \centering
        \includegraphics[width=0.55\linewidth]{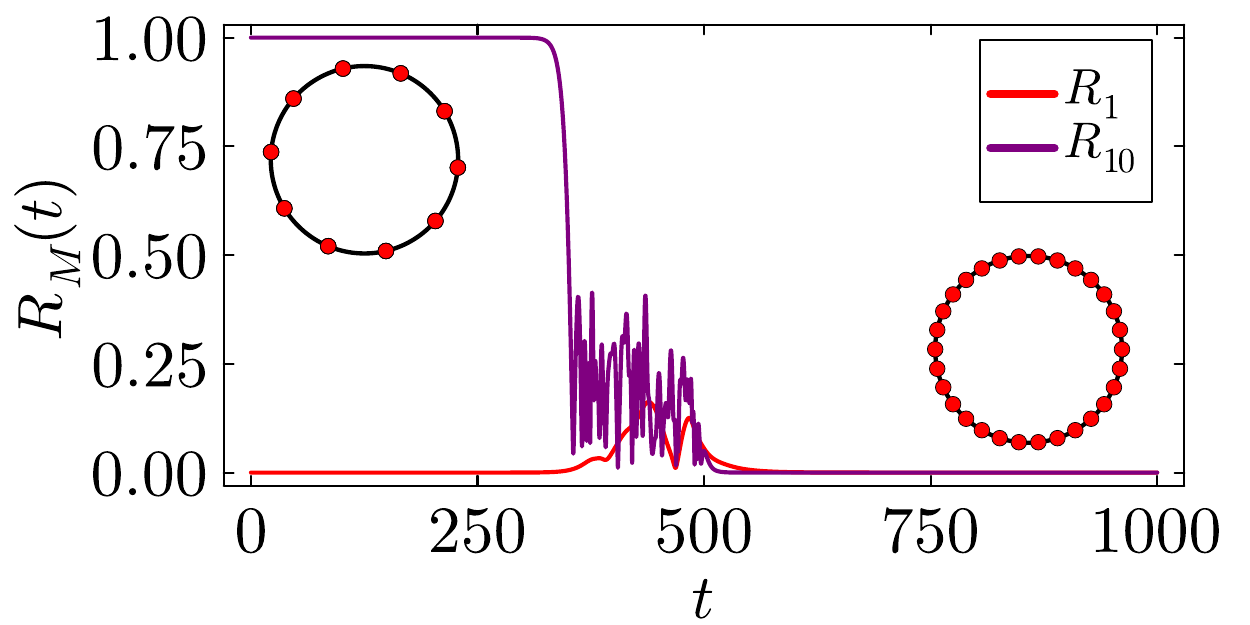}
        \caption{$v^{-1} = 0.6$}
        \label{subfig:twisted state emergence a}
    \end{subfigure}

    \medskip

    \begin{subfigure}[t]{\columnwidth}
        \centering
        \includegraphics[width=0.55\linewidth]{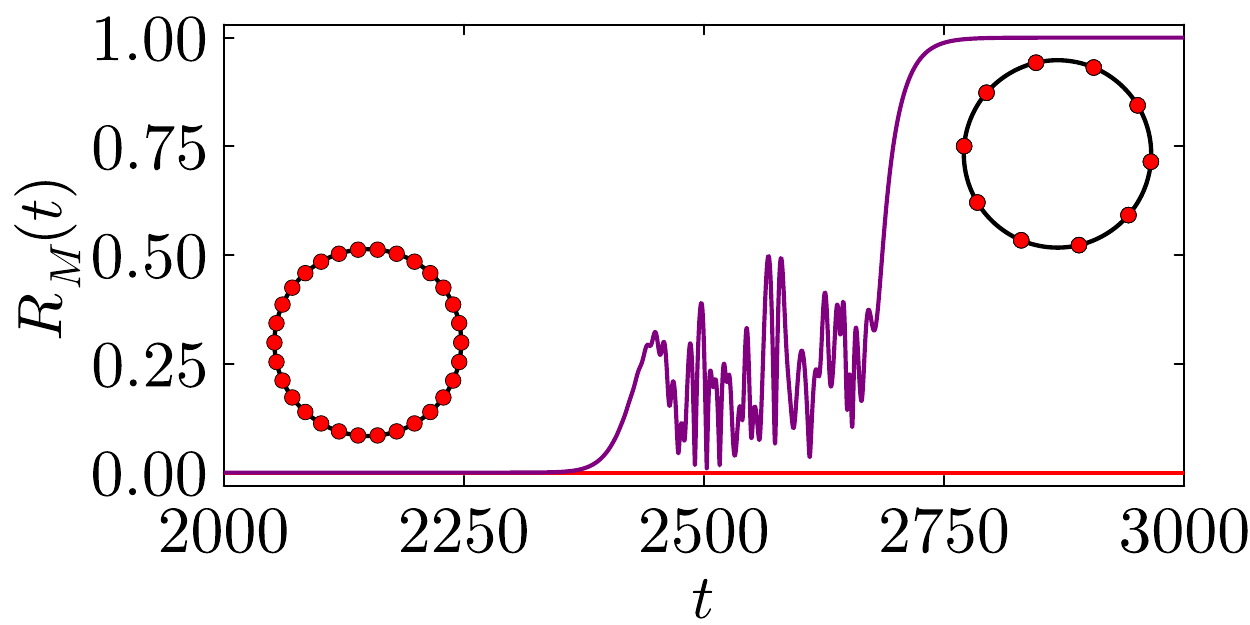}
        \caption{$v^{-1} = 1.0$}
        \label{subfig:twisted state emergence b}
    \end{subfigure}

    \caption{ Evolution of the Kuramoto--Daido order parameters $R_1$ and $R_{10}$ at two values of $v^{-1}$ illustrating transitions between a splay state and a $3$-twisted state. (a) Emergence of a stable splay state from an unstable $3$-twisted state. Initial condition $\theta_i(0) = 2\pi \cdot 3(i-1) / 30$, for $i = 1,...,30$ and history function $h_i(t) = 0.92t + \theta_i(0)$. (b) Emergence of a stable $3$-twisted state from an unstable splay state. Initial condition $\theta_i(0) = 2\pi (i-1)  / 30 $ and history function $h_i(t) = 0.95t + \theta_i(0)$. Parameters and connectivity function $w(x)$ as in Fig.\ref{fig:twisted branches}.}
    
    \label{fig:twisted state emergence}
\end{figure}

We consider an example of a $3$-twisted state and the splay ($1$-twisted) state in a network with $N=30$, since the $3$-twisted state is stable for a large range of $v^{-1}$ and seen in simulations. The solution branches for small delays are shown in Fig.~\ref{fig:twisted branches}. There are parameter values where the splay state is stable while the $3$-twisted state unstable (e.g. $v^{-1} = 0.6$), and others where the $3$-twisted state is stable but the splay state is unstable (e.g. $v^{-1} = 1.0$). The $3$-twisted state in this network consists of ten equally spaced phase clusters, each containing three oscillators. A natural set of statistics to describe these states are the Kuramoto-Daido order parameters $R_M = \left| \frac{1}{N} \sum_{n=1}^N e^{i M \theta_n(t)} \right|$. A characteristic property of twisted states is that $R_1 = 0$ for all states with $q \neq 0$. More generally, if a solution consists of $M$ equally spaced clusters, then $R_M = 1$. Since the $3$-twisted state forms ten clusters, we therefore track the evolution of $R_1(t)$ and $R_{10}(t)$ to identify changes between solutions. An example is shown in Fig.~\ref{fig:twisted state emergence}. For $v^{-1} = 0.6$, the system is initialized near the unstable $3$-twisted state shown in Fig.~\ref{fig:twisted branches} using a history function with the phase configuration of the $3$-twisted state but the frequency of the stable splay state ($\Omega =0.92$). The system evolves to the stable splay state. For $v^{-1} = 1.0$, the system is initialized near the unstable splay state, with frequency of the $3$-twisted state ($\Omega =0.95$), and it evolves to the stable $3$-twisted state.

\subsection{Patterns of phase-locked states}

The symmetry properties of the ring network can be used to identify additional phase-locked solutions. The ring network of $N$ nodes has dihedral symmetry $D_N$, which allows us to apply the $H/K$ theorem of equivariant bifurcation theory \cite{Golubitsky2002}  to find (up to conjugacy) all possible $2\pi$-periodic solution patterns which the network may support. For example, with $N=6$ the possible $1:1$ phase-locked solutions are of the form shown in Table \ref{tab:ring N=6 patterns}.  A brief description of the methodology illustrating the main arguments from symmetric bifurcation theory used to produce Table \ref{tab:ring N=6 patterns} is provided in Appendix \ref{appendixsymmetry}.  Note that solutions corresponding to other subgroup pairs $(H,K)$ are possible, but these are $1:2$ phase-locked states and not listed in Table \ref{tab:ring N=6 patterns}.

\begin{figure}[t]
    \centering
    \begin{subfigure}[t]{0.2\textwidth}
        \centering
        \tikzset{every picture/.style={line width=0.75pt}} 

\begin{tikzpicture}[x=0.75pt,y=0.75pt,yscale=-1,xscale=1]

\draw   (96.56,139) .. controls (96.56,103.61) and (125.42,74.93) .. (161.02,74.93) .. controls (196.62,74.93) and (225.48,103.61) .. (225.48,139) .. controls (225.48,174.39) and (196.62,203.07) .. (161.02,203.07) .. controls (125.42,203.07) and (96.56,174.39) .. (96.56,139) -- cycle ;
\draw  [color={rgb, 255:red, 0; green, 0; blue, 255 }  ,draw opacity=1 ][fill={rgb, 255:red, 255; green, 255; blue, 255 }  ,fill opacity=1 ] (139.97,74.93) .. controls (139.97,63.37) and (149.39,54) .. (161.02,54) .. controls (172.65,54) and (182.08,63.37) .. (182.08,74.93) .. controls (182.08,86.49) and (172.65,95.86) .. (161.02,95.86) .. controls (149.39,95.86) and (139.97,86.49) .. (139.97,74.93) -- cycle ;
\draw  [color={rgb, 255:red, 0; green, 0; blue, 255 }  ,draw opacity=1 ][fill={rgb, 255:red, 255; green, 255; blue, 255 }  ,fill opacity=1 ] (139.97,203.07) .. controls (139.97,191.51) and (149.39,182.14) .. (161.02,182.14) .. controls (172.65,182.14) and (182.08,191.51) .. (182.08,203.07) .. controls (182.08,214.63) and (172.65,224) .. (161.02,224) .. controls (149.39,224) and (139.97,214.63) .. (139.97,203.07) -- cycle ;
\draw  [color={rgb, 255:red, 0; green, 0; blue, 255 }  ,draw opacity=1 ][fill={rgb, 255:red, 255; green, 255; blue, 255 }  ,fill opacity=1 ] (84.37,106.41) .. controls (84.37,94.85) and (93.8,85.48) .. (105.43,85.48) .. controls (117.05,85.48) and (126.48,94.85) .. (126.48,106.41) .. controls (126.48,117.97) and (117.05,127.34) .. (105.43,127.34) .. controls (93.8,127.34) and (84.37,117.97) .. (84.37,106.41) -- cycle ;
\draw  [color={rgb, 255:red, 0; green, 0; blue, 255 }  ,draw opacity=1 ][fill={rgb, 255:red, 255; green, 255; blue, 255 }  ,fill opacity=1 ] (195.2,107.05) .. controls (195.2,95.49) and (204.62,86.12) .. (216.25,86.12) .. controls (227.88,86.12) and (237.31,95.49) .. (237.31,107.05) .. controls (237.31,118.61) and (227.88,127.98) .. (216.25,127.98) .. controls (204.62,127.98) and (195.2,118.61) .. (195.2,107.05) -- cycle ;
\draw  [color={rgb, 255:red, 0; green, 0; blue, 255 }  ,draw opacity=1 ][fill={rgb, 255:red, 255; green, 255; blue, 255 }  ,fill opacity=1 ] (84,170.01) .. controls (84,158.45) and (93.43,149.08) .. (105.06,149.08) .. controls (116.69,149.08) and (126.11,158.45) .. (126.11,170.01) .. controls (126.11,181.57) and (116.69,190.94) .. (105.06,190.94) .. controls (93.43,190.94) and (84,181.57) .. (84,170.01) -- cycle ;
\draw  [color={rgb, 255:red, 0; green, 0; blue, 255 }  ,draw opacity=1 ][fill={rgb, 255:red, 255; green, 255; blue, 255 }  ,fill opacity=1 ] (194.83,170.65) .. controls (194.83,159.09) and (204.25,149.72) .. (215.88,149.72) .. controls (227.51,149.72) and (236.94,159.09) .. (236.94,170.65) .. controls (236.94,182.21) and (227.51,191.58) .. (215.88,191.58) .. controls (204.25,191.58) and (194.83,182.21) .. (194.83,170.65) -- cycle ;
\draw  [color={rgb, 255:red, 0; green, 0; blue, 0 }  ,draw opacity=1 ][fill={rgb, 255:red, 1; green, 0; blue, 0 }  ,fill opacity=1 ] (155.61,54) .. controls (155.61,51.01) and (158.03,48.58) .. (161.02,48.58) .. controls (164.01,48.58) and (166.44,51.01) .. (166.44,54) .. controls (166.44,56.99) and (164.01,59.42) .. (161.02,59.42) .. controls (158.03,59.42) and (155.61,56.99) .. (155.61,54) -- cycle ;
\draw    (161.02,54) -- (161.02,74.93) ;
\draw  [color={rgb, 255:red, 0; green, 0; blue, 0 }  ,draw opacity=1 ][fill={rgb, 255:red, 1; green, 0; blue, 0 }  ,fill opacity=1 ] (100.01,85.48) .. controls (100.01,82.49) and (102.43,80.07) .. (105.43,80.07) .. controls (108.42,80.07) and (110.84,82.49) .. (110.84,85.48) .. controls (110.84,88.47) and (108.42,90.9) .. (105.43,90.9) .. controls (102.43,90.9) and (100.01,88.47) .. (100.01,85.48) -- cycle ;
\draw  [color={rgb, 255:red, 0; green, 0; blue, 0 }  ,draw opacity=1 ][fill={rgb, 255:red, 1; green, 0; blue, 0 }  ,fill opacity=1 ] (210.84,86.12) .. controls (210.84,83.13) and (213.26,80.7) .. (216.25,80.7) .. controls (219.24,80.7) and (221.67,83.13) .. (221.67,86.12) .. controls (221.67,89.11) and (219.24,91.53) .. (216.25,91.53) .. controls (213.26,91.53) and (210.84,89.11) .. (210.84,86.12) -- cycle ;
\draw  [color={rgb, 255:red, 0; green, 0; blue, 0 }  ,draw opacity=1 ][fill={rgb, 255:red, 1; green, 0; blue, 0 }  ,fill opacity=1 ] (99.64,149.08) .. controls (99.64,146.09) and (102.07,143.67) .. (105.06,143.67) .. controls (108.05,143.67) and (110.47,146.09) .. (110.47,149.08) .. controls (110.47,152.07) and (108.05,154.5) .. (105.06,154.5) .. controls (102.07,154.5) and (99.64,152.07) .. (99.64,149.08) -- cycle ;
\draw  [draw opacity=0] (155.61,224) .. controls (155.61,221.01) and (158.03,218.58) .. (161.02,218.58) .. controls (164.01,218.58) and (166.44,221.01) .. (166.44,224) .. controls (166.44,226.99) and (164.01,229.42) .. (161.02,229.42) .. controls (158.03,229.42) and (155.61,226.99) .. (155.61,224) -- cycle ;
\draw  [color={rgb, 255:red, 0; green, 0; blue, 0 }  ,draw opacity=1 ][fill={rgb, 255:red, 1; green, 0; blue, 0 }  ,fill opacity=1 ] (210.47,149.72) .. controls (210.47,146.72) and (212.89,144.3) .. (215.88,144.3) .. controls (218.87,144.3) and (221.3,146.72) .. (221.3,149.72) .. controls (221.3,152.71) and (218.87,155.13) .. (215.88,155.13) .. controls (212.89,155.13) and (210.47,152.71) .. (210.47,149.72) -- cycle ;
\draw    (105.43,85.48) -- (105.43,106.41) ;
\draw    (216.25,86.12) -- (216.25,107.05) ;
\draw    (105.06,149.08) -- (105.06,170.01) ;
\draw    (161.02,182.14) -- (161.02,203.07) ;
\draw    (215.88,149.72) -- (215.88,170.65) ;
\draw  [draw opacity=0] (142.96,43.14) .. controls (148.22,39.68) and (154.4,37.68) .. (161.02,37.68) .. controls (167.86,37.68) and (174.24,39.81) .. (179.61,43.49) -- (161.02,74.93) -- cycle ; \draw   (142.96,43.14) .. controls (148.22,39.68) and (154.4,37.68) .. (161.02,37.68) .. controls (167.86,37.68) and (174.24,39.81) .. (179.61,43.49) ;  
\draw    (142.96,43.14) -- (140.86,44.43) ;
\draw [shift={(138.31,46)}, rotate = 328.42] [fill={rgb, 255:red, 0; green, 0; blue, 0 }  ][line width=0.08]  [draw opacity=0] (10.72,-5.15) -- (0,0) -- (10.72,5.15) -- (7.12,0) -- cycle    ;
\draw  [color={rgb, 255:red, 0; green, 0; blue, 0 }  ,draw opacity=1 ][fill={rgb, 255:red, 1; green, 0; blue, 0 }  ,fill opacity=1 ] (155.61,182.14) .. controls (155.61,179.15) and (158.03,176.72) .. (161.02,176.72) .. controls (164.01,176.72) and (166.44,179.15) .. (166.44,182.14) .. controls (166.44,185.13) and (164.01,187.56) .. (161.02,187.56) .. controls (158.03,187.56) and (155.61,185.13) .. (155.61,182.14) -- cycle ;

\draw (150,10) node [font=\huge][anchor=north west][inner sep=0.75pt]   [align=left] {$ \Omega $};

\end{tikzpicture}
        \caption{Pattern 1}
        \label{fig:a}
    \end{subfigure}
    \hspace{0.01\textwidth}
    \begin{subfigure}[t]{0.2\textwidth}
        \centering
        \tikzset{every picture/.style={line width=0.75pt}} 

\begin{tikzpicture}[x=0.75pt,y=0.75pt,yscale=-1,xscale=1]

\draw   (126.56,134) .. controls (126.56,98.61) and (155.42,69.93) .. (191.02,69.93) .. controls (226.62,69.93) and (255.48,98.61) .. (255.48,134) .. controls (255.48,169.39) and (226.62,198.07) .. (191.02,198.07) .. controls (155.42,198.07) and (126.56,169.39) .. (126.56,134) -- cycle ;
\draw  [color={rgb, 255:red, 0; green, 0; blue, 255 }  ,draw opacity=1 ][fill={rgb, 255:red, 255; green, 255; blue, 255 }  ,fill opacity=1 ] (169.97,69.93) .. controls (169.97,58.37) and (179.39,49) .. (191.02,49) .. controls (202.65,49) and (212.08,58.37) .. (212.08,69.93) .. controls (212.08,81.49) and (202.65,90.86) .. (191.02,90.86) .. controls (179.39,90.86) and (169.97,81.49) .. (169.97,69.93) -- cycle ;
\draw  [color={rgb, 255:red, 0; green, 0; blue, 255 }  ,draw opacity=1 ][fill={rgb, 255:red, 255; green, 255; blue, 255 }  ,fill opacity=1 ] (169.97,198.07) .. controls (169.97,186.51) and (179.39,177.14) .. (191.02,177.14) .. controls (202.65,177.14) and (212.08,186.51) .. (212.08,198.07) .. controls (212.08,209.63) and (202.65,219) .. (191.02,219) .. controls (179.39,219) and (169.97,209.63) .. (169.97,198.07) -- cycle ;
\draw  [color={rgb, 255:red, 0; green, 0; blue, 255 }  ,draw opacity=1 ][fill={rgb, 255:red, 255; green, 255; blue, 255 }  ,fill opacity=1 ] (114.37,101.41) .. controls (114.37,89.85) and (123.8,80.48) .. (135.43,80.48) .. controls (147.05,80.48) and (156.48,89.85) .. (156.48,101.41) .. controls (156.48,112.97) and (147.05,122.34) .. (135.43,122.34) .. controls (123.8,122.34) and (114.37,112.97) .. (114.37,101.41) -- cycle ;
\draw  [color={rgb, 255:red, 0; green, 0; blue, 255 }  ,draw opacity=1 ][fill={rgb, 255:red, 255; green, 255; blue, 255 }  ,fill opacity=1 ] (225.2,102.05) .. controls (225.2,90.49) and (234.62,81.12) .. (246.25,81.12) .. controls (257.88,81.12) and (267.31,90.49) .. (267.31,102.05) .. controls (267.31,113.61) and (257.88,122.98) .. (246.25,122.98) .. controls (234.62,122.98) and (225.2,113.61) .. (225.2,102.05) -- cycle ;
\draw  [color={rgb, 255:red, 0; green, 0; blue, 255 }  ,draw opacity=1 ][fill={rgb, 255:red, 255; green, 255; blue, 255 }  ,fill opacity=1 ] (114,165.01) .. controls (114,153.45) and (123.43,144.08) .. (135.06,144.08) .. controls (146.69,144.08) and (156.11,153.45) .. (156.11,165.01) .. controls (156.11,176.57) and (146.69,185.94) .. (135.06,185.94) .. controls (123.43,185.94) and (114,176.57) .. (114,165.01) -- cycle ;
\draw  [color={rgb, 255:red, 0; green, 0; blue, 255 }  ,draw opacity=1 ][fill={rgb, 255:red, 255; green, 255; blue, 255 }  ,fill opacity=1 ] (224.83,165.65) .. controls (224.83,154.09) and (234.25,144.72) .. (245.88,144.72) .. controls (257.51,144.72) and (266.94,154.09) .. (266.94,165.65) .. controls (266.94,177.21) and (257.51,186.58) .. (245.88,186.58) .. controls (234.25,186.58) and (224.83,177.21) .. (224.83,165.65) -- cycle ;
\draw  [color={rgb, 255:red, 0; green, 0; blue, 0 }  ,draw opacity=1 ][fill={rgb, 255:red, 1; green, 0; blue, 0 }  ,fill opacity=1 ] (185.61,49) .. controls (185.61,46.01) and (188.03,43.58) .. (191.02,43.58) .. controls (194.01,43.58) and (196.44,46.01) .. (196.44,49) .. controls (196.44,51.99) and (194.01,54.42) .. (191.02,54.42) .. controls (188.03,54.42) and (185.61,51.99) .. (185.61,49) -- cycle ;
\draw    (191.02,49) -- (191.02,69.93) ;
\draw  [color={rgb, 255:red, 0; green, 0; blue, 0 }  ,draw opacity=1 ][fill={rgb, 255:red, 1; green, 0; blue, 0 }  ,fill opacity=1 ] (130.01,122.34) .. controls (130.01,119.35) and (132.43,116.93) .. (135.43,116.93) .. controls (138.42,116.93) and (140.84,119.35) .. (140.84,122.34) .. controls (140.84,125.33) and (138.42,127.76) .. (135.43,127.76) .. controls (132.43,127.76) and (130.01,125.33) .. (130.01,122.34) -- cycle ;
\draw  [color={rgb, 255:red, 0; green, 0; blue, 0 }  ,draw opacity=1 ][fill={rgb, 255:red, 1; green, 0; blue, 0 }  ,fill opacity=1 ] (240.84,122.98) .. controls (240.84,119.98) and (243.26,117.56) .. (246.25,117.56) .. controls (249.24,117.56) and (251.67,119.98) .. (251.67,122.98) .. controls (251.67,125.97) and (249.24,128.39) .. (246.25,128.39) .. controls (243.26,128.39) and (240.84,125.97) .. (240.84,122.98) -- cycle ;
\draw  [color={rgb, 255:red, 0; green, 0; blue, 0 }  ,draw opacity=1 ][fill={rgb, 255:red, 1; green, 0; blue, 0 }  ,fill opacity=1 ] (129.64,144.08) .. controls (129.64,141.09) and (132.07,138.67) .. (135.06,138.67) .. controls (138.05,138.67) and (140.47,141.09) .. (140.47,144.08) .. controls (140.47,147.07) and (138.05,149.5) .. (135.06,149.5) .. controls (132.07,149.5) and (129.64,147.07) .. (129.64,144.08) -- cycle ;
\draw  [color={rgb, 255:red, 0; green, 0; blue, 0 }  ,draw opacity=1 ][fill={rgb, 255:red, 1; green, 0; blue, 0 }  ,fill opacity=1 ] (185.61,219) .. controls (185.61,216.01) and (188.03,213.58) .. (191.02,213.58) .. controls (194.01,213.58) and (196.44,216.01) .. (196.44,219) .. controls (196.44,221.99) and (194.01,224.42) .. (191.02,224.42) .. controls (188.03,224.42) and (185.61,221.99) .. (185.61,219) -- cycle ;
\draw  [color={rgb, 255:red, 0; green, 0; blue, 0 }  ,draw opacity=1 ][fill={rgb, 255:red, 1; green, 0; blue, 0 }  ,fill opacity=1 ] (240.47,144.72) .. controls (240.47,141.72) and (242.89,139.3) .. (245.88,139.3) .. controls (248.87,139.3) and (251.3,141.72) .. (251.3,144.72) .. controls (251.3,147.71) and (248.87,150.13) .. (245.88,150.13) .. controls (242.89,150.13) and (240.47,147.71) .. (240.47,144.72) -- cycle ;
\draw    (135.43,101.41) -- (135.43,122.34) ;
\draw    (246.25,102.05) -- (246.25,122.98) ;
\draw    (135.06,144.08) -- (135.06,165.01) ;
\draw    (191.02,198.07) -- (191.02,219) ;
\draw    (245.88,144.72) -- (245.88,165.65) ;
\draw  [draw opacity=0] (172.96,38.14) .. controls (178.22,34.68) and (184.4,32.68) .. (191.02,32.68) .. controls (197.86,32.68) and (204.24,34.81) .. (209.61,38.49) -- (191.02,69.93) -- cycle ; \draw   (172.96,38.14) .. controls (178.22,34.68) and (184.4,32.68) .. (191.02,32.68) .. controls (197.86,32.68) and (204.24,34.81) .. (209.61,38.49) ;  
\draw    (172.96,38.14) -- (170.86,39.43) ;
\draw [shift={(168.31,41)}, rotate = 328.42] [fill={rgb, 255:red, 0; green, 0; blue, 0 }  ][line width=0.08]  [draw opacity=0] (10.72,-5.15) -- (0,0) -- (10.72,5.15) -- (7.12,0) -- cycle    ;

\draw (180,5) node [font=\huge][anchor=north west][inner sep=0.75pt]   [align=left] {$ \Omega $};

\end{tikzpicture}
        \caption{Pattern 2}
        \label{fig:b}
    \end{subfigure}
    \hspace{0.01\textwidth}
    \begin{subfigure}[t]{0.2\textwidth}
        \centering
        \tikzset{every picture/.style={line width=0.75pt}} 

\begin{tikzpicture}[x=0.75pt,y=0.75pt,yscale=-1,xscale=1]

\draw   (265.25,125) .. controls (265.25,89.61) and (294.11,60.93) .. (329.71,60.93) .. controls (365.31,60.93) and (394.17,89.61) .. (394.17,125) .. controls (394.17,160.39) and (365.31,189.07) .. (329.71,189.07) .. controls (294.11,189.07) and (265.25,160.39) .. (265.25,125) -- cycle ;
\draw  [color={rgb, 255:red, 0; green, 0; blue, 255 }  ,draw opacity=1 ][fill={rgb, 255:red, 255; green, 255; blue, 255 }  ,fill opacity=1 ] (308.66,60.93) .. controls (308.66,49.37) and (318.08,40) .. (329.71,40) .. controls (341.34,40) and (350.77,49.37) .. (350.77,60.93) .. controls (350.77,72.49) and (341.34,81.86) .. (329.71,81.86) .. controls (318.08,81.86) and (308.66,72.49) .. (308.66,60.93) -- cycle ;
\draw  [color={rgb, 255:red, 0; green, 0; blue, 255 }  ,draw opacity=1 ][fill={rgb, 255:red, 255; green, 255; blue, 255 }  ,fill opacity=1 ] (308.66,189.07) .. controls (308.66,177.51) and (318.08,168.14) .. (329.71,168.14) .. controls (341.34,168.14) and (350.77,177.51) .. (350.77,189.07) .. controls (350.77,200.63) and (341.34,210) .. (329.71,210) .. controls (318.08,210) and (308.66,200.63) .. (308.66,189.07) -- cycle ;
\draw  [color={rgb, 255:red, 0; green, 0; blue, 255 }  ,draw opacity=1 ][fill={rgb, 255:red, 255; green, 255; blue, 255 }  ,fill opacity=1 ] (253.06,92.41) .. controls (253.06,80.85) and (262.49,71.48) .. (274.12,71.48) .. controls (285.75,71.48) and (295.17,80.85) .. (295.17,92.41) .. controls (295.17,103.97) and (285.75,113.34) .. (274.12,113.34) .. controls (262.49,113.34) and (253.06,103.97) .. (253.06,92.41) -- cycle ;
\draw  [color={rgb, 255:red, 0; green, 0; blue, 255 }  ,draw opacity=1 ][fill={rgb, 255:red, 255; green, 255; blue, 255 }  ,fill opacity=1 ] (363.89,93.05) .. controls (363.89,81.49) and (373.31,72.12) .. (384.94,72.12) .. controls (396.57,72.12) and (406,81.49) .. (406,93.05) .. controls (406,104.61) and (396.57,113.98) .. (384.94,113.98) .. controls (373.31,113.98) and (363.89,104.61) .. (363.89,93.05) -- cycle ;
\draw  [color={rgb, 255:red, 0; green, 0; blue, 255 }  ,draw opacity=1 ][fill={rgb, 255:red, 255; green, 255; blue, 255 }  ,fill opacity=1 ] (252.69,156.01) .. controls (252.69,144.45) and (262.12,135.08) .. (273.75,135.08) .. controls (285.38,135.08) and (294.8,144.45) .. (294.8,156.01) .. controls (294.8,167.57) and (285.38,176.94) .. (273.75,176.94) .. controls (262.12,176.94) and (252.69,167.57) .. (252.69,156.01) -- cycle ;
\draw  [color={rgb, 255:red, 0; green, 0; blue, 255 }  ,draw opacity=1 ][fill={rgb, 255:red, 255; green, 255; blue, 255 }  ,fill opacity=1 ] (363.52,156.65) .. controls (363.52,145.09) and (372.95,135.72) .. (384.57,135.72) .. controls (396.2,135.72) and (405.63,145.09) .. (405.63,156.65) .. controls (405.63,168.21) and (396.2,177.58) .. (384.57,177.58) .. controls (372.95,177.58) and (363.52,168.21) .. (363.52,156.65) -- cycle ;
\draw  [color={rgb, 255:red, 0; green, 0; blue, 0 }  ,draw opacity=1 ][fill={rgb, 255:red, 1; green, 0; blue, 0 }  ,fill opacity=1 ] (324.3,40) .. controls (324.3,37.01) and (326.72,34.58) .. (329.71,34.58) .. controls (332.71,34.58) and (335.13,37.01) .. (335.13,40) .. controls (335.13,42.99) and (332.71,45.42) .. (329.71,45.42) .. controls (326.72,45.42) and (324.3,42.99) .. (324.3,40) -- cycle ;
\draw    (329.71,40) -- (329.71,60.93) ;
\draw  [color={rgb, 255:red, 0; green, 0; blue, 0 }  ,draw opacity=1 ][fill={rgb, 255:red, 1; green, 0; blue, 0 }  ,fill opacity=1 ] (251.7,104.48) .. controls (251.7,101.49) and (254.13,99.07) .. (257.12,99.07) .. controls (260.11,99.07) and (262.53,101.49) .. (262.53,104.48) .. controls (262.53,107.47) and (260.11,109.9) .. (257.12,109.9) .. controls (254.13,109.9) and (251.7,107.47) .. (251.7,104.48) -- cycle ;
\draw  [color={rgb, 255:red, 0; green, 0; blue, 0 }  ,draw opacity=1 ][fill={rgb, 255:red, 1; green, 0; blue, 0 }  ,fill opacity=1 ] (268.33,135.08) .. controls (268.33,132.09) and (270.76,129.67) .. (273.75,129.67) .. controls (276.74,129.67) and (279.16,132.09) .. (279.16,135.08) .. controls (279.16,138.07) and (276.74,140.5) .. (273.75,140.5) .. controls (270.76,140.5) and (268.33,138.07) .. (268.33,135.08) -- cycle ;
\draw  [color={rgb, 255:red, 0; green, 0; blue, 0 }  ,draw opacity=1 ][fill={rgb, 255:red, 1; green, 0; blue, 0 }  ,fill opacity=1 ] (379.16,135.72) .. controls (379.16,132.72) and (381.58,130.3) .. (384.57,130.3) .. controls (387.57,130.3) and (389.99,132.72) .. (389.99,135.72) .. controls (389.99,138.71) and (387.57,141.13) .. (384.57,141.13) .. controls (381.58,141.13) and (379.16,138.71) .. (379.16,135.72) -- cycle ;
\draw    (274.12,92.41) -- (257.12,104.48) ;
\draw    (273.75,135.08) -- (273.75,156.01) ;
\draw    (384.57,135.72) -- (384.57,156.65) ;
\draw  [draw opacity=0] (311.65,29.14) .. controls (316.91,25.68) and (323.09,23.68) .. (329.71,23.68) .. controls (336.55,23.68) and (342.93,25.81) .. (348.3,29.49) -- (329.71,60.93) -- cycle ; \draw   (311.65,29.14) .. controls (316.91,25.68) and (323.09,23.68) .. (329.71,23.68) .. controls (336.55,23.68) and (342.93,25.81) .. (348.3,29.49) ;  
\draw    (311.65,29.14) -- (309.56,30.43) ;
\draw [shift={(307,32)}, rotate = 328.42] [fill={rgb, 255:red, 0; green, 0; blue, 0 }  ][line width=0.08]  [draw opacity=0] (10.72,-5.15) -- (0,0) -- (10.72,5.15) -- (7.12,0) -- cycle    ;
\draw    (329.71,189.07) -- (312.71,201.14) ;
\draw    (384.94,93.05) -- (367.94,105.12) ;
\draw  [color={rgb, 255:red, 0; green, 0; blue, 0 }  ,draw opacity=1 ][fill={rgb, 255:red, 1; green, 0; blue, 0 }  ,fill opacity=1 ] (362.53,105.12) .. controls (362.53,102.13) and (364.95,99.7) .. (367.94,99.7) .. controls (370.93,99.7) and (373.36,102.13) .. (373.36,105.12) .. controls (373.36,108.11) and (370.93,110.53) .. (367.94,110.53) .. controls (364.95,110.53) and (362.53,108.11) .. (362.53,105.12) -- cycle ;
\draw  [color={rgb, 255:red, 0; green, 0; blue, 0 }  ,draw opacity=1 ][fill={rgb, 255:red, 1; green, 0; blue, 0 }  ,fill opacity=1 ] (307.3,201.14) .. controls (307.3,198.15) and (309.72,195.72) .. (312.71,195.72) .. controls (315.71,195.72) and (318.13,198.15) .. (318.13,201.14) .. controls (318.13,204.13) and (315.71,206.56) .. (312.71,206.56) .. controls (309.72,206.56) and (307.3,204.13) .. (307.3,201.14) -- cycle ;
\draw  [dash pattern={on 2.25pt off 2.25pt}]  (274.12,92.41) -- (274,55) ;
\draw  [draw opacity=0] (241.4,101.16) .. controls (240.6,98.38) and (240.17,95.44) .. (240.17,92.41) .. controls (240.17,76.35) and (252.19,62.99) .. (268.04,60.23) -- (274.12,92.41) -- cycle ; \draw   (241.4,101.16) .. controls (240.6,98.38) and (240.17,95.44) .. (240.17,92.41) .. controls (240.17,76.35) and (252.19,62.99) .. (268.04,60.23) ;  
\draw    (241.4,101.16) -- (242.78,104.26) ;
\draw [shift={(244,107)}, rotate = 246.02] [fill={rgb, 255:red, 0; green, 0; blue, 0 }  ][line width=0.08]  [draw opacity=0] (10.72,-5.15) -- (0,0) -- (10.72,5.15) -- (7.12,0) -- cycle    ;

\draw  [draw opacity=0] (324.3,210) .. controls (324.3,207.01) and (326.72,204.58) .. (329.71,204.58) .. controls (332.71,204.58) and (335.13,207.01) .. (335.13,210) .. controls (335.13,212.99) and (332.71,215.42) .. (329.71,215.42) .. controls (326.72,215.42) and (324.3,212.99) .. (324.3,210) -- cycle ;

\draw (318.69,-4) node [font=\huge][anchor=north west][inner sep=0.75pt]   [align=left] {$ \Omega $};
\draw (222,57) node [font=\huge][anchor=north west][inner sep=0.75pt]   [align=left] {$ \psi $};

\end{tikzpicture}
        \caption{Pattern 3}
        \label{fig:c}
    \end{subfigure}

    \caption{Diagrams illustrating Patterns 1, 2, and 3 in Table \ref{tab:ring N=6 patterns}.
    For Pattern~3 in (c), \( y(t) = x(t + \psi) \) where $\psi$ is a phase lag which is determined by the particulars of the governing equations i.e., the interaction function $H$, weight function $w$ and delays depending on $v^{-1}$. }
    \label{fig:table patterns}
\end{figure}
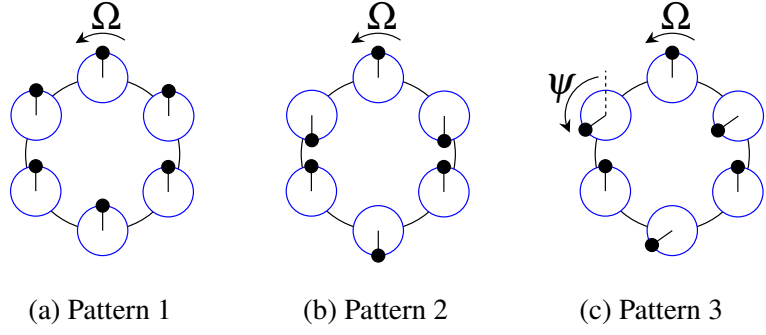

 \begin{widetext}

\begin{table}[h!]
    \centering
\begin{tabular}{|c|c|c|l|}
\hline
\# & \textbf{H} & \textbf{K} & \textbf{Pattern} \\
\hline
1  & $D_6$ & $D_6$ & $(x(t), x(t), x(t), x(t), x(t), x(t))$ \\
\hline
2  & $D_6$ & $S^1_3$ & $(x(t), x(t+\pi), x(t), x(t+\pi), x(t), x(t+\pi))$ \\
3  & $S^1_3$ & $S^1_3$ & $(x(t), y(t), x(t), y(t), x(t), y(t))$ \\
\hline
4  & $\mathbb{Z}_2 \times \mathbb{Z}_2$ & $\mathbb{Z}_2 \times \mathbb{Z}_2$ & $(x(t), x(t), y(t), x(t), x(t), y(t))$ \\
\hline
5  & $\mathbb{Z}_2(\kappa\rho )$ & $\mathbb{Z}_2(\kappa \rho )$ & $(x(t), x(t), y(t), z(t), z(t), y(t))$ \\
\hline
6  & $\mathbb{Z}_2 \times \mathbb{Z}_2$ & $\mathbb{Z}_2(\kappa)$ & $(x(t), y(t), y(t+\pi), x(t+\pi), y(t+\pi), y(t))$ \\
7  & $\mathbb{Z}_2(\kappa)$ & $\mathbb{Z}_2(\kappa)$ & $(x(t), y(t), z(t), u(t), z(t), y(t))$ \\
\hline
8 & $\mathbb{Z}_6$ & $\mathbb{Z}_2(\rho^3)$ & $(x(t), x(t+\tfrac{2\pi}{3}), x(t+\tfrac{4\pi}{3}), x(t), x(t+\tfrac{2\pi}{3}), x(t+\tfrac{4\pi}{3}))$ \\
9 & $\mathbb{Z}_2(\rho^3)$ & $\mathbb{Z}_2(\rho^3)$ & $(x(t), y(t), z(t), x(t), y(t), z(t))$ \\
\hline
10 & $\mathbb{Z}_6$ & $\mathbf{1}$ & $(x(t), x(t+\tfrac{2\pi}{6}), x(t+\tfrac{4\pi}{6}), x(t+\tfrac{6\pi}{6}), x(t+\tfrac{8\pi}{6}), x(t+\tfrac{10\pi}{6}))$ \\
11 & $\mathbb{Z}_3$ & $\mathbf{1}$ & $(x(t), y(t), x(t+\tfrac{2\pi}{3}), y(t+\tfrac{2\pi}{3}), x(t+\tfrac{4\pi}{3}), y(t+\tfrac{4\pi}{3}))$ \\
12 & $\mathbb{Z}_2(\kappa \rho)$ & $\mathbf{1}$ & $(x(t), x(t+\pi), y(t), z(t), z(t+\pi), y(t+\pi))$ \\
13 & $\mathbb{Z}_2(\rho^3)$ & $\mathbf{1}$ & $(x(t), y(t), z(t), x(t+\pi), y(t+\pi), z(t+\pi))$ \\
14 & $\mathbf{1}$ & $\mathbf{1}$ & $(x(t), y(t), z(t), u(t), v(t), w(t))$ \\
\hline
\end{tabular}
    \caption{The $1:1$ phase-locked solution patterns which may exist by symmetry for a six node ring network. $H$ and $K$ are subgroups of $D_6$ satisfying Lemma \ref{lem:HK}. Each pair defines a conjugacy class of solution patterns which may exist within the ring network depending on particulars of the model. We give an example pattern for each pair. See Appendix \ref{appendixsymmetry} for details on how these are computed.  }
    \label{tab:ring N=6 patterns}
\end{table}

\end{widetext}

\begin{figure}
\centering
\includegraphics[width=0.65\linewidth]{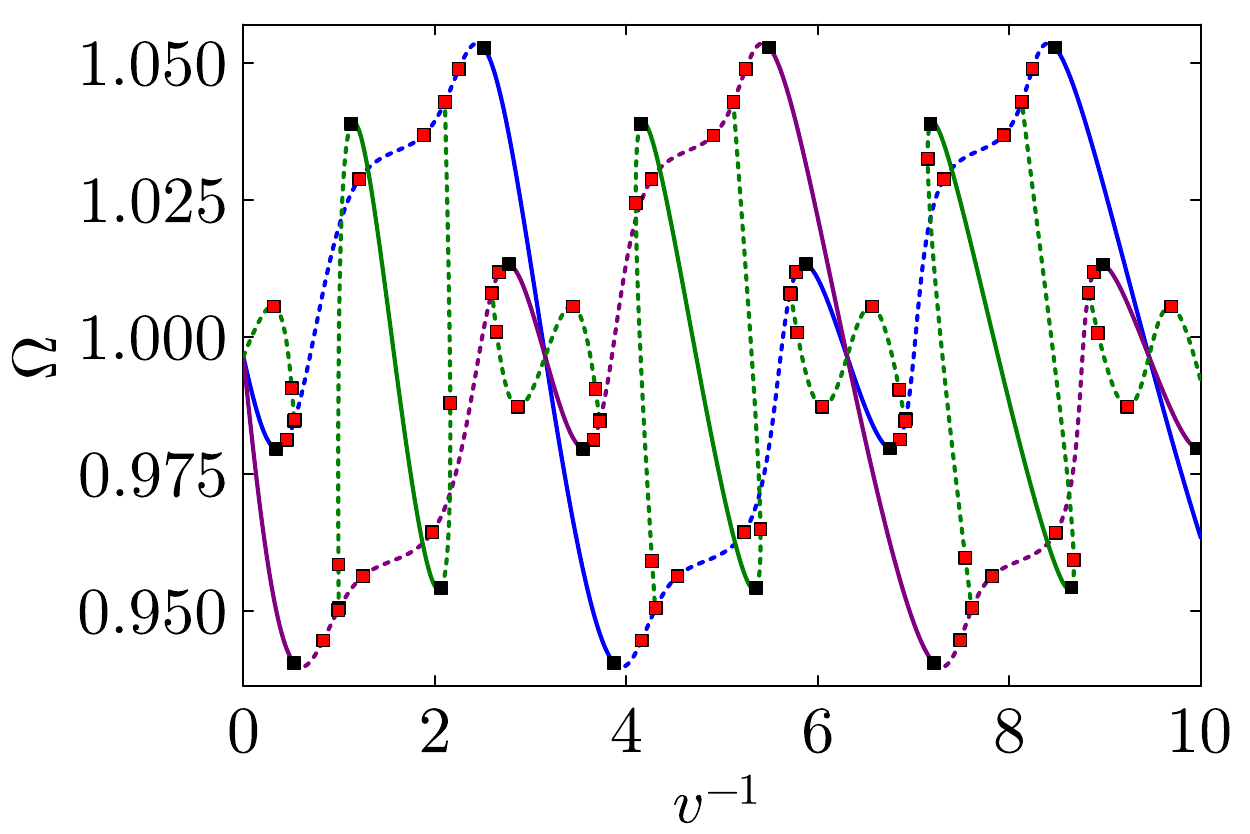}
\caption{Emergent frequency $\Omega$ as a function of  $v^{-1}$ in a six node ring network with distance dependent interactions. Three branches of phase-locked solutions from Table \ref{tab:ring N=6 patterns} are shown, but others also exist (not shown for clarity). 
Solutions for pattern 1 are in blue, pattern 2 in purple and pattern 3 in green. 
Stable solutions are indicated by solid lines and unstable solutions are dashed lines. 
Square markers denote bifurcations where a real eigenvalue crosses zero: black squares indicate changes of branch stability and red squares show other bifurcations where an eigenvalue crossing occurs for an already unstable solution. Some of these may be symmetry breaking bifurcations to additional solution branches from Table \ref{tab:ring N=6 patterns} which are not shown. 
Here $w(x) = \exp\left( -|x| \right)/2$.
Parameters:  $\omega=1$, $a=0$, $r=1$, $\sigma=0.1$, $d=1$ and $\tau_0=1$.
\label{Fig: table branches}
}
\end{figure}

Some of the solution types in Table \ref{tab:ring N=6 patterns} are phase locked states we have already discussed: Pattern 1 is synchrony, Pattern 2 is a $3$-twisted state, Pattern 8 is a $2$-twisted state and Pattern 10 is the splay state. 
Diagrams illustrating solution classes 1, 2 and 3 are shown in Fig.~\ref{fig:table patterns}. Solution class 3 corresponds to pattern 1 when $\psi=0$ and to pattern 2 when $\psi=\pi$. In Fig.~\ref{Fig: table branches} we show bifurcation diagrams in $v^{-1}$ for these three solution types for a particular choice of the other parameters.  The branches of solutions associated with pattern 3 emerge from bifurcations of the solution curves of either pattern 1 or pattern 2. 
Observe that there are values of $v^{-1}$ for which none of the plotted solutions are stable, for example at $v^{-1} = 1$. Here, numerical simulations suggest that stable solution branches include phase-locked solutions corresponding to Patterns 5,6,7,8 and 10.


\section{White matter plasticity\label{plasticity}}

It is now known that white matter is \textit{plastic}, namely that the myelination of an axon is activity dependent 
see e.g.,  \cite{Fields2008,Gibson2014,Vivo2019}.  Thus, given that myelination affects conduction speed, the delays that we have considered as fixed up until now may themselves evolve in a state dependent fashion.
The time-scale for white matter plasticity can be as short as days to weeks in the context of a new learning experience, though ranges up to years for e.g., changes associated with learning to read \cite{Huber2018}.
This has recently begun to be explored from a modeling perspective by Lefebvre and colleagues \cite{Noori2020,Park2020,Talidou2021,Talidou2022,Lefebvre2025}, with recent work by Fields and colleagues emphasizing the role of oligodendrocyte-mediated myelin plasticity in facilitating neural synchronization, whereby 
white matter plasticity may help distant brain regions remain synchronized by compensating for long transmission times through increased conduction speeds \cite{Pajevic2014,Pajevic2023}. 

Following the work of Lefebvre and colleagues it is natural to describe white matter plasticity in terms of a biologically motivated rule for the evolution of conduction speed.  For a fixed axonal fiber length between nodes $i$ and $j$ of length $d_{ij}$, we will consider a uniformly myelinated system that gives rise to a delay $\tau_{ij}(t) = d_{ij}/v_{ij}(t)$ for a signal traveling along the fiber with speed $v_{ij}$.  In the absence of any myelination we will define this speed to be $v^0_{ij}$, such that with an increase in myelination $v_{ij} > v_{ij}^0$.  To develop a plasticity rule it is important to bear in mind that i) myelination is activity dependent (via axon-glia interaction), and ii) axonal plasticity is slow and has a metabolic cost.  Within the phase oscillator context presented here we propose a state-dependent DDE model of the form
\begin{equation}
\frac{1}{\alpha}\FD{}{t} v_{ij}(t) = v_{ij}^0 - v_{ij}(t) + \kappa \int_{t-\tau_{ij}(t)}^t \mathcal{F}(\theta_j(s), \dot{\theta}_j(s)) \d s .
\label{WMPlasticity}
\end{equation}
Here, $\alpha^{-1}$ is a time-scale that can set the slowness of the plasticity with respect to the period of a single node oscillation.  The first term on the right hand side of (\ref{WMPlasticity}) sets the baseline for speed, namely that of the unmyelinated axon. The final term has an integral over the time that it takes for a signal to propagate along a fiber.  This is meant to help capture the fact that metabolic cost will depend upon activity over the whole length of the fiber whilst it is active.  
The strength of this activity is captured by the function $\mathcal{F}(\theta, \dot{\theta})$, and since this is expected to be strongly dependent on the firing frequency of the transmitting node we shall take $\mathcal{F}(\theta, \dot{\theta}) = \dot{\theta}$.  In this case, the plasticity rule (\ref{WMPlasticity}) takes the form
\begin{equation}
\frac{1}{\alpha}\FD{}{t} v_{ij}(t) = v_{ij}^0 - v_{ij}(t) + \kappa \left [ \theta_j(t) - \theta_j(t-d_{ij}/v_{ij}(t)) \right ].
\label{WMPlasticity1}
\end{equation}
In contrast to previous white matter plasticity rules for phase-oscillator networks as in \cite{Karimian2019,Park2020} the evolution for speed is governed by a state-dependent DDE, as opposed to an ODE. Of course, it is also possible to consider more traditional plasticity rules for the strength of connections, as in Hebbian learning, though this is not our focus here, although see \cite{Timms2014} for such approaches in coupled Kuramoto networks with delays.

\subsection{Phase-locked states revisited}

If the speed equation \eqref{WMPlasticity1} settles to a stationary time-independent value $\overline{v}_{ij}$, so that time-delays become fixed to $\overline{\tau}_{ij} = d_{ij}/\overline{v}_{ij}$, then we may expect to see the emergence of phase-locked states such as those discussed in \S \ref{phaselocked}.
Assuming solutions of the form $\theta_{i}(t) = \Omega t + \phi_i$, this will be the case if there is a solution to the set of equations
\begin{align}
\overline{v}_{ij} &= {v}_{ij}^0 + \kappa \Omega d_{ij}/\overline{v}_{ij},  \\
\Omega & = \omega + \sigma \sum_{j=1}^N W_{ij} H \left (\phi_j - \phi_i -\Omega d_{ij}/\overline{v}_{ij} \right ) ,
\end{align}
for $i,j = 1,\ldots,N$.
Fixing a phase, say $\phi_1=0$, gives a system of $N^2+N$ equations for the $N^2+N$ unknowns $\Omega, \phi_{i\neq1}, \overline{v}_{ij}$.  Assuming a solution with $\overline{v}_{ij} - v_{ij}^0 \geq 0$ exists,
we introduce small perturbations to the phase-locked state and write $\theta_i(t) = \Omega t + \phi_i + \delta \theta_i(t)$, and $v_{ij}(t) = \overline{v}_{ij}(t) +\delta v_{ij}(t)$, with corresponding perturbations 
$\tau_{ij}(t) = \overline{\tau}_{ij}(t) +\delta \tau_{ij}(t)$ with $\delta \tau_{ij} = - d_{ij}/(\overline{v}_{ij}
)^2 \delta v_{ij}$. Substitution into (\ref{PhaseNetwork}) and (\ref{WMPlasticity1}) and expanding to first order we obtain a set of linearized equations with solutions of the form $(\delta \theta_i(t), \delta \tau_{ij}(t)) = \e^{\lambda t} (u_i, r_{ij})$ for some non-zero amplitudes $u_i$ and $r_{ij}$ and $\lambda \in \CSet$.  These satisfy the linear equations
\begin{align}
&- \frac{(\overline{v}_{ij})^2}{d_{ij}} \left ( \frac{\lambda}{\alpha} + 1 \right ) r_{ij} = \kappa \left ( \Omega r_{ij}+u_j  ( 1-\e^{-\lambda \overline{\tau}_{ij}} ) \right ) ,\label{u}\\
&\lambda u_i = \sum_{j=1}^N  \left [ \mathcal{J}_{ij}(\lambda) u_{j} - \mathcal{J}_{ij}(0) (u_i+\Omega r_{ij})  \right ] ,\label{r}
\end{align}
where $\mathcal{J}_{ij}(\lambda) = \sigma W_{ij} H'(\phi_j-\phi_i -\Omega \overline{\tau}_{ij}) \e^{-\lambda \overline{\tau}_{ij}}$.  Rearranging (\ref{u}) for $r_{ij}$ allows us to eliminate it from (\ref{r}) to obtain
$\lambda  u_i = \sum_j \mathcal{B}_{ij} (\lambda) u_j$ where
\begin{equation}
\mathcal{B}_{ij}(\lambda) =  \mathcal{J}_{ij}(\lambda) - \delta_{ij} k_i + \frac{\Omega \kappa  \mathcal{J}_{ij}(0) ( 1-\e^{-\lambda \overline{\tau}_{ij}} )}{ \frac{(\overline{v}_{ij})^2}{d_{ij}}\left ( \frac{\lambda}{\alpha} + 1 \right ) + \Omega \kappa} ,
\end{equation}
and $k_i = \sum_k  \mathcal{J}_{ik}(0)$.  
For non-trivial solutions we must have that $\mathcal{E}(\lambda) = 0$, where
$\mathcal{E}(\lambda) = \det (\lambda I_N - \mathcal{B}(\lambda) )$.
We note that $\lambda=0$ is a always a solution of $\mathcal{E}(\lambda) = 0$, reflecting the invariance of solutions to phase-shifts.  In this way we may determine the stability of a phase locked state according to the condition $\text{Re}\, (\lambda) <0 $ (excluding the zero eigenvalue that arises from phase-shift symmetry).

\begin{figure}
\centering
\includegraphics[width=0.55\linewidth]{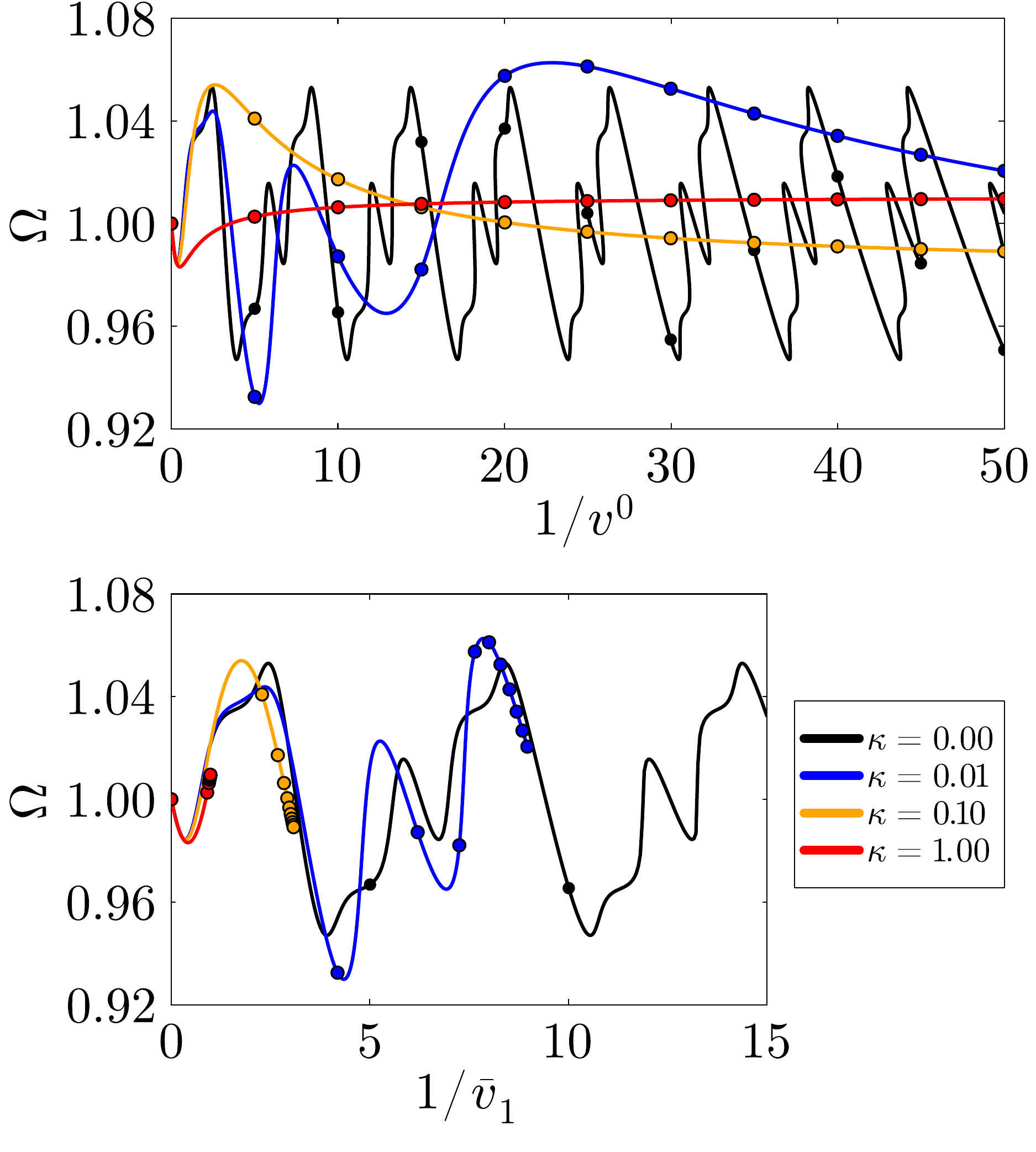}
\caption{Emergent frequency $\Omega$ as a function of $1/v^0$ and $1/\overline{v}_1$ for synchronous solutions in a six node ring network with distant dependent interactions and white matter plasticity at varying plasticity strengths.
Solutions in both plots were evaluated for $1/v^0 \in [0, 50]$. Circle markers denote solutions evaluated at uniformly spaced values of $1/v^0$. 
Parameters as in Fig.~\ref{Fig: table branches}.
\label{Fig:WMRing}
}
\end{figure}

By way of illustration of the effects of the plasticity rule (\ref{WMPlasticity1}) on the properties of phase-locked states we first revisit the architecture studied in \S \ref{ring}, namely a ring network, and set $v_{ij}^0 = v^0$ for all $(i,j)$. In this network, the conduction speed matrix $\overline{v}_{ij}$ inherits the circulant structure of the ring network with $\overline{v}_{ij} = \overline{v}_{\operatorname{dist} (i,j)}$ and rows generated by $\overline{v}_{j} = \overline{v}_{\operatorname{dist} (0,j)}$ for $j=0,...,5$. Branches of synchronous solutions for various values of the strength of white matter plasticity $\kappa$ are shown in Fig.~\ref{Fig:WMRing}. The top panel shows emergent frequency as a function of the parameter $1/v^0$. Here we see that plasticity promotes a \textit{flattening} of solution branches in the sense that emergent solutions all tend to have a roughly similar frequency for different values of $v^0$ and that the multiplicity of solution collective frequencies is suppressed in favor of just a single one for each phase-locked state.  The bottom panel shows emergent frequency as a function of $1/\overline{v}_{1}$, a quantity proportional to a resultant delay in this network. As the plasticity strength $\kappa$ increases solutions exist only for progressively smaller values of $1/\overline{v}_{1}$. The circle markers cluster near the end points of each branch, indicating that the system converges to the same solution independently of the initial speed $v^0$. This effect strengths with increasing plasticity strength.
Other phase-locked states show similar behavior to that shown in these two panels.

To establish that this phenomenon is not just specific to a ring network we next turn to more realistic brain-like architectures built from human connectome data.

\subsection{A connectome study}

\begin{figure}
\centering
\includegraphics[width=0.7\linewidth]{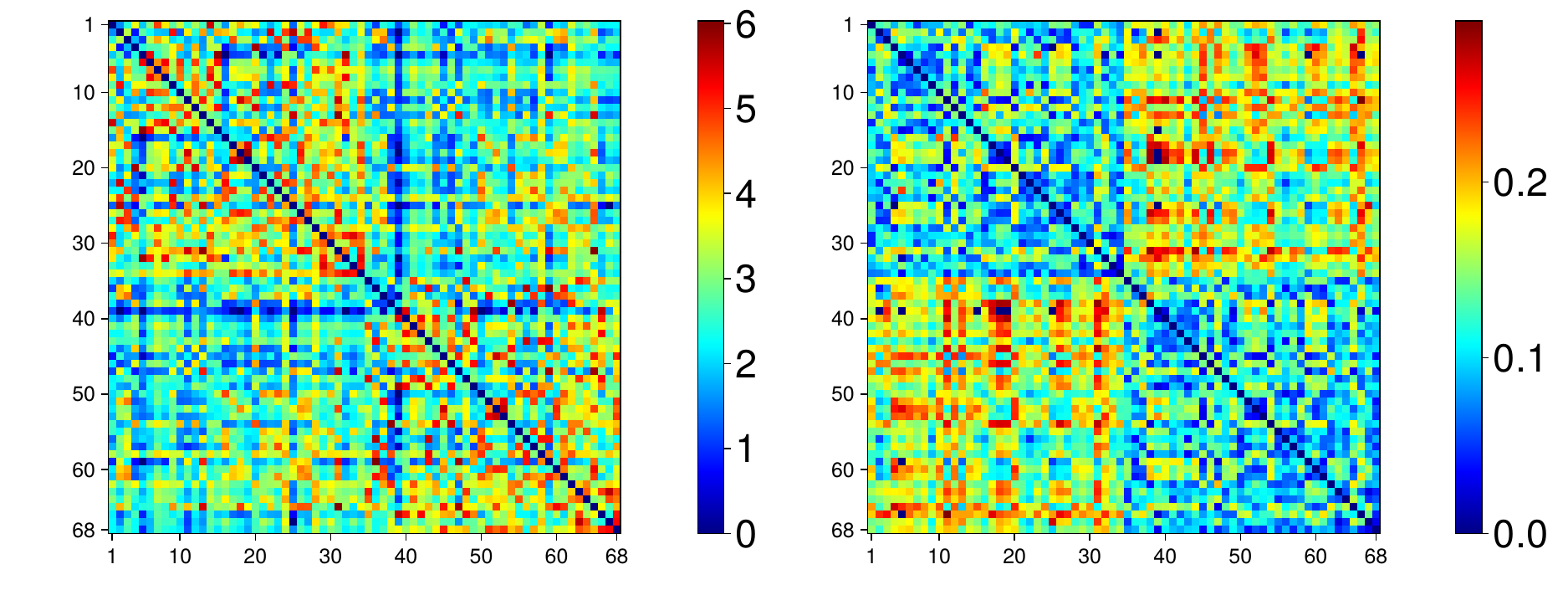}
\caption{Data from the Human Connectome Project, downsampled onto a $68$ node network. Left: Structural connectivity between nodes. Right: Axonal distance between nodes.}

\label{Fig:Connectome}

\end{figure}

Here we make use of space-time connectome data from the Human Connectome Project \cite{VanEssen2013}, downsampled onto a 68 node network using the Desikan--Killiany atlas as in Ref.~\onlinecite{Forrester2024}.  This gives us a symmetric set of positive structural connections $W_{ij}$ and a set of axonal distances $d_{ij}$ for long-range excitatory connections.  The weighted graphs of these connections are shown in Fig.~\ref{Fig:Connectome} using a matrix representation.

From this connectome data we can build a large scale brain model according to (\ref{PhaseNetwork}) with $\tau_{ij}=d_{ij}/v_{ij}$ and then explore the consequences of the white matter plasticity rule (\ref{WMPlasticity1}) on emergent network dynamics.  Here we shall make use of direct numerical simulations and not restrict ourselves to the study of simple phase-locked states.  To quantify network behavior it is useful to consider a measure of \textit{functional connectivity} (FC)
that can capture patterns of correlation and coherence between nodes based on temporal similarity.  A common choice for this, and one which we shall adopt is the pairwise phase-locking value (taken at an instant in time):
\begin{equation}\label{eq:plv}
F_{ij} = \frac{1}{2} \left ( 1 + \cos(\theta_j-\theta_i) \right ) .
\end{equation}

\begin{figure}
\centering
\includegraphics[width=0.7\linewidth]{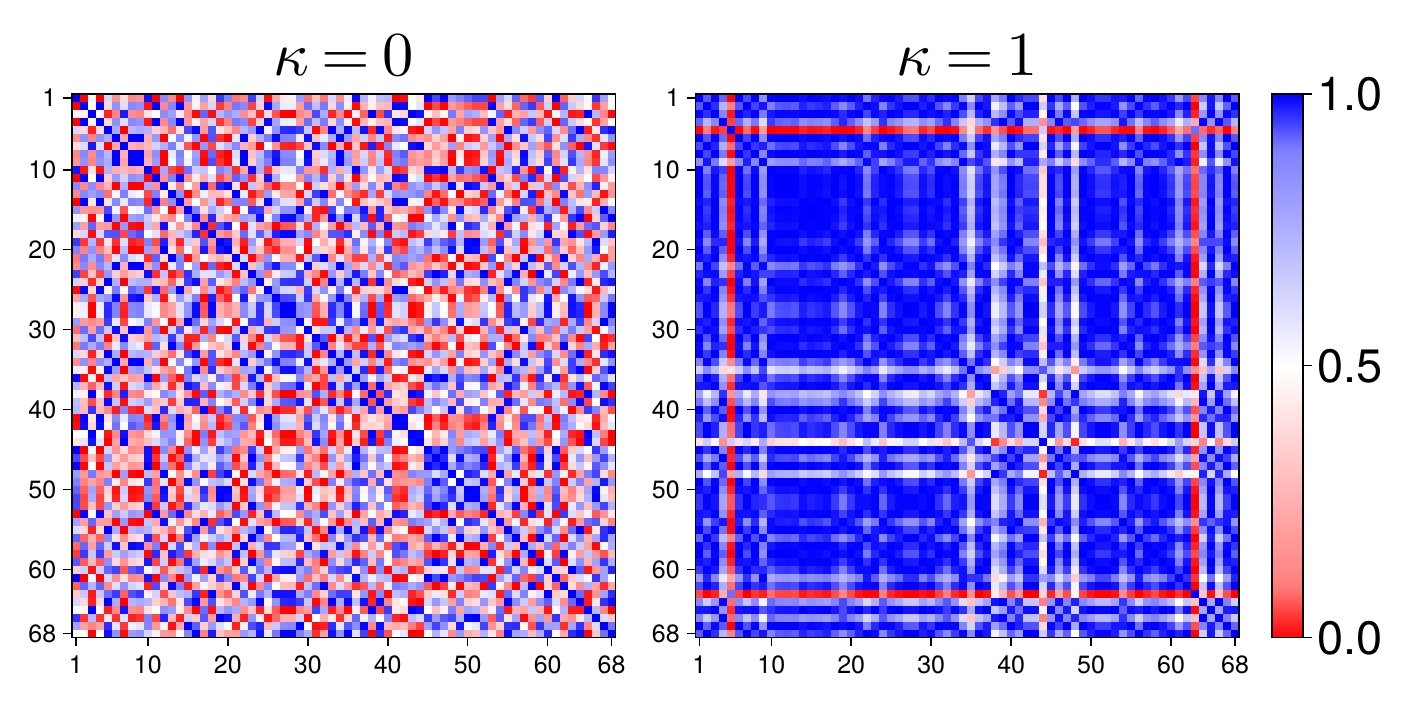}
\caption{Patterns of FC before ($\kappa=0$) and after ($\kappa=1$) the application of the white matter plasticity rule in an $N=68$ network. The color map shows values of $F_{ij}$ (\ref{eq:plv}) for each pair of nodes, blue indicates synchronous behavior and red indicates anti-synchronous behavior. The connectivity weights $W_{ij}$ and distances $d_{ij}$ are obtained from connectome data, see Figure \ref{Fig:Connectome}.
Initial conditions are given by $\theta_i(0) = 2\pi(i-1)/N$ and $v_{ij}(0)=v^0$, with $v^0=1.0$. History functions are $h_i(t) = t + \theta_i(0)$ for the phases and $h_{ij}(t) = v_{ij}(0)$ for the speeds. Simulations were run up to $t=1000$. 
The phase interaction function is $H(\theta) = K_1\left(a_0 -\sin(\theta - a_1) +r\sin(2(\theta-a_2))\right)$ with $K_1 = 6.778127193,\; a_0 = 0.780097488,\; a_1 = -0.640668464,\; a_2 = 0.770292062,\; r = 0.175334914$.
Other parameters: $\omega=1$, $\sigma =0.1$, $\alpha=0.1$.}

\label{Fig:FC}

\end{figure}

In Fig.~\ref{Fig:FC} we show a typical plot of FC with fixed delays and how the pattern of FC changes when these are allowed to evolve according to the white matter plasticity rule.  For the former we set $v_{ij}=v^0$ for all $(i,j)$ and for the latter we use this as initial data for $v_{ij}(0)$.

For the system with fixed delays ($\kappa = 0$), the FC matrix is largely unstructured. The values of $F_{ij}$ span almost the full range from $0$ to $1$, with many intermediate values, showing a wide range of pairwise correlations.  There are no clusters or groups of strongly correlated nodes, suggesting the network does not exhibit synchronized activity. 

However, after the application of the white matter plasticity rule ($\kappa=1$), the FC matrix has a much stronger correlation structure. Many pairs of nodes satisfy $F_{ij}=1$, indicating  that these nodes are exactly synchronized. Only a small number of nodes are anti-synchronous with the rest of the network, while the remaining values of $F_{ij}$ are close to $1$. This connectome data study shows that the plasticity rule can promote the emergence of a highly synchronized network.

Overall, both the ring network and the connectome studies show that white matter plasticity can result in more organized collective behavior. In the ring network, plasticity reduced multistability and dependence of the end dynamics on the initial conductance speeds, while the connectome study showed it can increase synchronization. Together, these demonstrate that the white matter plasticity rule can simplify network dynamics by driving the system to fewer and more coherent states.

\section{Discussion\label{discussion}}

Networks of coupled oscillators with delayed interactions are often described by systems of DDEs.  The analysis of phase-locked states in such models is challenging, because (i) limit cycles are rarely available in closed form, (ii) studying their linear stability requires Floquet theory for delay equations \cite{sieber2011}.  However, when reducing such systems to phase-oscillator networks, which is expected to hold when interactions are weak, the resulting set of DDEs is far more tractable.  This is because oscillatory states in the un-reduced model correspond to relative equilibria in the reduced model.  Although this is still a DDE system, the construction and linear stability of solutions is readily performed. We have developed the mathematical framework for this in the first part of this paper and illustrated its effectiveness with a number of bifurcation studies of various networks.  The second part of this paper recognizes that in a neural context there is a growing interest in the effects of white matter plasticity on the dynamics of brain networks.  This naturally leads us to models of state-dependent DDEs, and the consideration of models whereby the delays between nodes can be modulated by the activity of the nodes.  Motivated by existing biological knowledge of white matter plasticity and previous modeling, we have presented a new candidate model of such a process.  Moreover, we have shown that the analysis of phase-locked states can be developed by building on our earlier framework.  Our modeling, analysis, and simulation studies lend credence to the idea that specific classes of white matter plasticity rules can help adjust network timing to favor more coherent synchronous states.

A number of possible immediate extensions of the work presented in this paper can also be considered.  Firstly, we have assumed a simple bi-harmonic form for the phase interaction function for illustrative purposes, whereas this can be derived from a biophysical model as a more general periodic function (as outlined in Appendix \ref{phasereduction}).  However, it is important to note that the theory we have developed here is general and applies to an arbitrary phase interaction function.  In the context of large scale brain modeling it would thus be of interest to revisit previous (more numerical) studies of neural mass networks, such as those in Ref.~\onlinecite{Forrester2024}, with the techniques presented here.  

Secondly, we can consider an improvement over the phase reduction using a \textit{phase-amplitude} reduction that retains some representation of distance from cycle.  There is a growing activity in this field, as exemplified by \cite{Castejon2013,Wilson2016,Monga2019,Ermentrout2019,Coombes2024,Nicks2024a}, and more recently it has been considered within the context of delay-induced oscillations \cite{Nicks2024}. In this framework, the multistability of phase-locked states of the same type but with different collective frequencies can be resolved as the existence of multiple distinct collective orbits with different oscillation amplitudes. Phase-amplitude reductions can therefore provide a route to understanding how delayed interactions shape dynamics near and away from a limit cycle, including off-cycle perturbations that are important for transitions between network states.

Thirdly, the approach presented here for networks can be generalized to cover continuum networks of the type recently discussed in 
\cite[Ch. 6.9]{Coombes2023}, and one might consider adapting the approach the approach in \cite{Kim2023} (for global coupling) to analyze spiral chimeras.  Also, it is well to mention that we have only considered identical oscillators and it would be instructive to introduce some heterogeneity and treat networks of non-identical oscillators.

More generally, the model of adaptive activity dependent myelination explored in this paper is one of an emerging class of DDE models for white matter plasticity. While our model supports the evolution of axonal speeds to stationary values, similar state dependent DDE rules with the inclusion of a sigmoidal function of the integrated activity over the fiber are investigated in \cite{Ruschel2026, Coombes2026} and have been found to generate network attractors in the form of explosive relaxation-type oscillations. The phenomenon can be understood through the lens of a slow-fast separation of time scales: Conduction speeds slowly evolve and network dynamics drift along stable branches of the bifurcation diagram for fixed speeds which acts as a frozen fast subsystem. Trajectories can approach slow fixed points, but can also switch between branches or undergo relaxation-type oscillations when the slow flow reverses direction between switches. Our work, along that of \cite{Ruschel2026, Coombes2026}, indicates that white matter plasticity can be utilized to reshape the attractor landscape of an oscillator network with the phase-locked states providing the network states and dynamic delays inducing changes in stability and transitions between them. This opens up the possibility that modification of communication delays could be a route to sculpting network attractors between oscillatory states in addition to the established route of changing connection weights \cite{Kori2001, Ashwin2024}. A realistic target when considering heterogeneous or connectome-like physical networks without symmetry is to sculpt excitable network attractors in which phase-locked states are stable but input-dependent transitions occur between them \cite{Ashwin2024}. Exploring the computational abilities of networks subject to delay plasticity presents an interesting avenue for future investigation.  

Finally, we highlight that the delay-induced effects on phase-locking introduced in this paper may also be useful for understanding the dynamical implications of degradation of myelin. Aging, as well as neurodegenerative and demyelinating diseases such as multiple sclerosis, leads to changes and damage in white matter, ultimately altering conduction speeds which is expected to disrupt coordinated neural communication \cite{sorrentino2022}. Developing white matter plasticity models for pathological myelin evolution and degradation may help us to understand disease progression in terms of the destabilization of coherent activity or changes in possible transitions between network states. Conversely, we may also be able to model the effects of remyelination therapies and the extent to which physiological timing relationships can be restored.

\begin{acknowledgments} This work was supported by The Leverhulme Trust through grant RPG-2025-052, ``White Matter Computation: Utilising axonal delays to sculpt network attractors'' 
and by the Scientific and Technological Research Council of Türkiye (TÜBİTAK) under the BIDEB 2219 International Postdoctoral Research Scholarship Program (Project No: 1059B192301947).
We would also like to thank J\'er\'emie Lefebvre for interesting discussions about white matter and its plasticity.
\end{acknowledgments}

\appendix

\section{Phase reduction\label{phasereduction}}

The theory of weakly coupled oscillators is a well established tool for the reduction of oscillatory networks, and see e.g., \cite[Ch. 6]{Coombes2023}.  Given that it is most often invoked in the absence of delays we briefly review the reduction process for a network with a space-time connectome as prescribed by (\ref{Network1}).

For $\sigma =0$ we introduce a phase on the cycle with dynamics $\dot{\theta}_i = \omega$, so that the dynamics can be reconstructed as $x_i(t) = s(\theta_i(t)/\omega)$  For sufficiently weak coupling the dynamics is expected to evolve as 
\begin{align}
&\FD{}{t} \theta_i(t) = \omega 
+ \sigma \sum_{j=1}^N W_{ij}  \left \langle Q(\theta_i(t)) G(s(\theta_i(t)/\omega), s(\theta_j(t -\tau_{ij})/\omega)) \right \rangle ,
\end{align}
where $Q \in \RSet^m$ is the $2 \pi$-periodic infinitesimal phase response curve, and $\langle \cdot \rangle$ denotes the standard inner product between vectors.
Moving to a co-rotating frame with the introduction of an angle $\phi_i(t) = \theta_i(t) - \omega t$, and recognizing that $\phi_i$ evolves slowly (namely $\dot{\phi}_i \sim O(\sigma)$), we may use first order averaging to turn the non-autonomous evolution equation for $\phi_i$ into an autonomous one given by
\begin{align}
&\FD{}{t} \phi_i (t)  \simeq \sigma \sum_{j=1}^N W_{ij} \frac{1}{T} \int_0^T \left \langle Q(\phi_i(t) +\omega t') \right. \times \nonumber \\
&\left. G(s(\phi_i(t)/\omega +t'), s(\phi_j(t -\tau_{ij})/\omega +t' - \tau_{ij})) \right \rangle \d t'.
\end{align}
In terms of the original phase variables this gives the dynamical system
\begin{align}
\FD{}{t} \theta_i (t) & = \omega + \sigma \sum_{j=1}^N W_{ij} H (\theta_j(t -\tau_{ij}) - \theta_i(t)) + O(\sigma^2),
\end{align}
where
\begin{equation}
H(\theta) = \frac{1}{2 \pi} \int_0^{2 \pi} \left \langle Q(t) G(s(t/\omega), s((\theta+t)/\omega)) \right \rangle \d t.
\label{H}
\end{equation}

\section{Parametric plotting\label{appendix parametric plotting}}

\subsection{The two oscillator network \label{appendix parametric plotting two node}}
The emergent frequency $\Omega$ for the two node network described in \S \ref{N=2} is given by 
\begin{align}
    \Omega &= \omega + \sigma H(\psi - \Omega \tau), \nonumber \\
    \Omega &= \omega + \sigma H( -\psi - \Omega \tau).
\label{eq:2node freq}
\end{align}
We plot the solution branches $\Omega = \Omega(\tau)$ and $\psi = \psi(\tau)$ by introducing the parametrization $\eta := \Omega \tau$. Subtracting the equations from (\ref{eq:2node freq}) leads to
\begin{equation}
    H(\psi - \eta) - H(-\psi - \eta) =0.
\label{eq: n2 plotting psi}
\end{equation}
For each value of \(\eta\), the phase difference \(\psi = \psi(\eta)\) is obtained numerically by solving \ref{eq: n2 plotting psi}. The emergent frequency $\Omega$ depends explicitly on $\eta$ since (\ref{eq:2node freq}) becomes 
$\Omega(\eta) = \omega + \sigma H(\psi(\eta) - \eta)$. We also rewrite $\tau$ in terms of $\eta$, 
\begin{equation}
    \tau(\eta) = \frac{\eta}{\Omega(\eta)} = \frac{\eta}{\omega + \sigma H(\psi(\eta) - \eta)},
\end{equation}
for $\Omega(\eta) \ne 0$. The curves $\Omega = \Omega(\tau)$ and $\psi = \psi(\tau)$ can therefore be plotted parametrically as 
\begin{equation}
    \left( \tau(\eta), \Omega(\eta)\right) = \left(\frac{\eta}{\omega + \sigma H(\psi(\eta) - \eta)} ,\omega + \sigma H(\psi(\eta) - \eta)\right), 
\label{parametric_2node}
\end{equation}
and 
$\left(\tau(\eta), \psi(\eta)\right)$
for $\eta \in\mathbb{R}$, where $\psi(\eta)$ is the solution to (\ref{eq: n2 plotting psi}). For each value of $\eta$, we numerically compute the eigenvalues of the solution using the bifurcation package DDE-BIFTOOL \cite{Engelborghs2002, Sieber2016}. These eigenvalues determine the linear stability of the solutions along the curve.

\subsection{Ring networks \label{appendix parametric plotting ring}}
The emergent frequency for several phase locked states in the ring network described in \S \ref{ring} can be expressed as 
\begin{equation}
    \Omega = \omega + \sigma \sum_{j=0}^{N-1} W_j\, H\!\left(\psi_j-\Omega \tau_j\right).
\label{eq:ring freq general}
\end{equation}
For synchrony $\psi_i =0$, for the splay state $\psi_i = 2\pi i/N$ and for a $q$-twisted state $\psi_i = 2\pi q i /N$.
We plot the curves parametrically by first introducing
\begin{equation}
    F(x) := x - \sigma W_0 H(-x\tau_0),
\end{equation}
so that (\ref{eq:ring freq general})  can be rewritten as
\begin{equation}
    F(\Omega) = \omega + \sigma \sum_{j=1}^{N-1} W_j\, H\!\left(\psi_j-\Omega \tau_j\right).
\end{equation}
Equivalently, 
\begin{equation}
    \Omega = F^{-1}\!\left(
    \omega + \sigma \sum_{j=1}^{N-1} W_j\, H\!\left(\psi_j-\Omega \tau_j\right)
    \right),
\end{equation}
provided that $F$ is invertible, with its inverse computed numerically. In all examples considered here, the parameters were chosen to ensure $F$ is invertible. 
We use the parametrization $\eta_1 := \Omega v^{-1}$ and since $\tau_j = \operatorname{dist} (0,j) v^{-1}$ for $j \ne 0$, we obtain
\begin{equation}
    \Omega(\eta_1) = F^{-1}\!\left(
    \omega + \sigma \sum_{j=1}^{N-1} W_j\, H\!\left(\psi_j- \operatorname{dist} (0,j) \eta_1 \right)
    \right).
\label{eq:plot splay}
\end{equation}
The solution curves $\Omega = \Omega(v^{-1})$ are plotted as $\left(\eta_1/\Omega(\eta_1), \Omega(\eta_1) \right)$, and the stability along these branches is computed using DDE-BIFTOOL.

\section{$H/K$ Theorem of equivariant bifurcation theory\label{appendixsymmetry}}

Consider a six node network with bidirectional coupling which can be circulant and/or delayed. The network has spatial symmetry group $D_6 = \langle \rho, \kappa \rangle$ where $ \rho = (1 2 3 4 5 6)$ represents a rotation and $\kappa = (12)(36)(45)$ represents a reflection. Let $X(t) = (x_1(t),...,x_6(t))$ is a $2\pi$-periodic solution of the system, $X(t) = X(t + 2\pi)$, satisfying the equations 
\begin{equation}
    \frac{\mathrm d}{ \mathrm{d}t} x_i = f(x_i) +  \sigma \sum_{j=1}^{N} W_{ij}  G \left( x_i(t),x_j(t-\tau_{ij}) \right),
\end{equation}
where each $x_i \in \mathbb{R}^m$ with $m\ge2$ ensuring the existence of a limit cycle to justify phase reduction to equation (\ref{PhaseNetwork}). 

A spatio-temporal symmetry of the solution $X(t)$ is defined by a pair $(\gamma, \theta) \in D_6 \times S^1$ such that $\gamma X(t + \theta) = X(t)$ for some $\gamma \in D_6$ and phase shift $\theta \in [0,2\pi]$. The full group of spatio-temporal symmetries of $X(t)$ forms a subgroup of $D_6 \times S^1$, which can be identified with a pair of subgroups $(H,K)$,  where $K \subset D_6$ is the spatial symmetry subgroup
\begin{equation}
    K = \{\gamma \in D_6 : \gamma X(t) = X(t) \quad \forall t\},
\end{equation}
and $H \subset D_6$ is the subgroup containing all spatial components of the spatio-temporal symmetries, defined by
\begin{equation}
    H = \{\gamma \in D_6 : \gamma \{X(t)\} = \{X(t)\}\}.
\end{equation}
There also exists a homeomorphism  $\Theta : H \to S^1$ whose kernel is $K$. By applying algebraic constraints on possible pairs 
$(H,K)$, we can classify the types of $2\pi$-periodic solutions that the network can support. In particular, Lemma 3.1 from \cite{Golubitsky2002} states  
\begin{lemma}  \label{lem:HK}
\begin{enumerate}
        \renewcommand{\labelenumi}{\alph{enumi}.} 
        \item $K$ is an isotropy subgroup for the action of $D_6$.
        \item $K$ is a normal subgroup of $H$ and $H/K$ is cyclic.
        \item If $\operatorname{dim} \operatorname{Fix}(K)= 2$, then either $H=K$ or $H=N(K)$, where $N(K)$ denotes the normalizer of $K$.
    \end{enumerate}
\end{lemma}
Using this Lemma, we find up to conjugacy all possible pairs $(H,K)$ for the $D_6$ action on the ring network and identify a corresponding example oscillation pattern. Table \ref{tab:ring N=6 patterns} summarizes these pairs and their associated $2\pi$-periodic solution patterns (restricted to patterns where all nodes share a common frequency).

\bibliography{ChaosPaperPhaseDelays}

\end{document}